\documentclass[12pt]{article}
\usepackage{amsmath}
\usepackage{amssymb}
\usepackage{vatola}

\def\q{\quad}
\def\qq{\qquad}
\def\mod{\pmod}
\def\t{\text}
\def\f{\frac}
\def\e{\equiv}
\def\b{\binom}
\def\a{\alpha}

\def\phq#1{\varphi(q^{#1})}
\def\psq#1{\psi(q^{#1})}
\def\qtq#1{\q\t{#1}\q}
\def\sls#1#2{(\f{#1}{#2})}
 \def\ls#1#2{\big(\f{#1}{#2}\big)}
\def\Ls#1#2{\Big(\f{#1}{#2}\Big)}
\let \pro=\proclaim
\let \endpro=\endproclaim

\begin{document}
 \par\q\newline
\centerline {\bf Ramanujan's theta functions and sums of
 triangular numbers}

$$\q$$
\centerline{Zhi-Hong Sun}
\par\q\newline
\centerline{School of Mathematical Sciences}
 \centerline{Huaiyin
Normal University} \centerline{Huaian,
Jiangsu 223300, P.R. China}
\centerline{Email: zhsun@hytc.edu.cn} \centerline{Homepage:
http://www.hytc.edu.cn/xsjl/szh}

 \abstract{Let $\Bbb Z$ and $\Bbb N$ be the set of integers
 and the set of positive integers, respectively. For
 $a_1,a_2,\ldots,a_k,n\in\Bbb N$ let $N(a_1,a_2,\ldots,a_k;n)$ be the number of
 representations of $n$ by $a_1x_1^2+a_2x_2^2+\cdots+a_kx_k^2$, and
 let $t(a_1,a_2,\ldots,a_k;n)$ be the number of
 representations of $n$ by $a_1\frac {x_1(x_1-1)}2+a_2\frac{x_2(x_2-1)}2+\cdots+a_k\f{x_k(x_k-1)}2$
  $(x_1,\ldots,x_k\in\Bbb Z$). In this paper, by using
 Ramanujan's theta functions $\varphi(q)$ and $\psi(q)$
 we reveal many relations
between $t(a_1,a_2,\ldots,a_k;n)$ and
$N(a_1,a_2,\ldots,a_k;8n+a_1+\cdots+a_k)$ for $k=3,4$.
 \par\q
 }
 \endabstract

\section*{1. Introduction}

\par Ramanujan's theta functions $\varphi(q)$ and $\psi(q)$ are defined
by
$$\varphi(q)=\sum_{n=-\infty}^{\infty}q^{n^2}=1+2\sum_{n=1}^{\infty}
q^{n^2}\qtq{and} \psi(q)=\sum_{n=0}^{\infty}q^{n(n+1)/2}\ (|q|<1).$$
From [BCH, Lemma 4.1] or [Be1] we know that for $|q|<1$,
 $$\align &\psi(q)^2=\varphi(q)\psi(q^2),\tag
 1.1\\&\varphi(q)=\varphi(q^4)+2q\psi(q^8),\tag 1.2
 \\&\varphi(q)^2=\phq 2^2+4q\psq 4^2=\phq 4^2+4q^2\psq 8^2+4q\psq 4^2,\tag 1.3
  \\&\psi(q)\psi(q^3)=\varphi(q^6)\psi(q^4)+q
 \varphi(q^2)\psi(q^{12}).\tag 1.4
 \endalign$$
By (1.2), for $k=1,2,3,\ldots$,
$$\varphi(q^k)=\varphi(q^{4k})+2q^k\psi(q^{8k})
=\varphi(q^{16k})+2q^{4k}\psi(q^{32k}) +2q^k\psi(q^{8k}).\tag 1.5$$

\par\q  Let $\Bbb Z$ and $\Bbb N$ be the set of integers
 and the set of positive integers, respectively.
  For $a_1,a_2,\ldots,
 a_k\in\Bbb N$ and $n\in\Bbb N \cup \{0\}$ set
$$\align &N(a_1,a_2,\ldots,a_k;n)=\big|\{(x_1,\ldots,x_k)\in \Bbb Z^k\ |
\ n=a_1x_1^2+a_2x_2^2+\cdots+a_kx_k^2 \}\big|,
\\&t(a_1,a_2,\ldots,a_k;n)\\&=\Big|\Big\{(x_1,\ldots,x_k)\in \Bbb
Z^k\ \big|\ n\ =a_1\f{x_1(x_1-1)}2+
a_2\f{x_2(x_2-1)}2+\cdots+a_k\f{x_k(x_k-1)}2\Big\}\Big|,
\\&t'(a_1,a_2,\ldots,a_k;n)\\&=\Big|\Big\{(x_1,\ldots,x_k)\in \Bbb
N^k\ \big|\ n\ =a_1\f{x_1(x_1-1)}2+
a_2\f{x_2(x_2-1)}2+\cdots+a_k\f{x_k(x_k-1)}2\Big\}\Big|.\endalign$$
The numbers $\f{x(x-1)}2\ (x\in\Bbb Z)$ are called triangular
numbers. Since $\f {x(x-1)}2=\f{(-x+1)(-x)}2$
we have
$$t(a_1,a_2,\ldots,a_k;n)=2^kt'(a_1,a_2,\ldots,
a_k;n).\tag 1.6$$ It is evident that for $|q|<1$,
$$\align&\sum_{n=0}^{\infty}N(a_1,\ldots,a_k;n)q^{n}=\varphi(q^{a_1})
\cdots\varphi(q^{a_k}),\tag 1.7
\\&\sum_{n=0}^{\infty}t'(a_1,\ldots,a_k;n)q^{n}=\psi(q^{a_1})
\cdots \psi(q^{a_k}).\tag 1.8\endalign$$
\par Let $a_1,\ldots,a_k\in\Bbb N$ and
$$C(a_1,\ldots,a_k)=2^k+2^{k-1}\Big(\f{i_1(i_1-1)(i_1-2)(i_1-3)}{4!}
+\f{i_1(i_1-1)i_2}2+i_1i_3\Big),$$ where $i_j$ denotes the number of
elements in $\{a_1,\ldots,a_k\}$ which are equal to $j$.
 In 2005 Adiga, Cooper and
Han [ACH] showed that
$$t'(a_1,a_2,\ldots,a_k;n)=\f {N(a_1,\ldots,a_k;8n+a_1+\cdots+a_k)}
{C(a_1,\ldots,a_k)}\q\t{for $a_1+\cdots+a_k\le 7$}.
 \tag 1.9$$ In
2008 Baruah, Cooper and Hirschhorn [BCH] proved that
$$\aligned &t'(a_1,a_2,\ldots,a_k;n)
\\&=\f {N(a_1,\ldots,a_k;8n+8)-N(a_1,\ldots,a_k;2n+2)}
{C(a_1,\ldots,a_k)}\q \t{for $a_1+\cdots+a_k=8$}.\endaligned \tag
1.10$$
\par For $n\in\Bbb N$ let
$$r_3(n)=N(1,1,1;n)=\big|\big\{(x,y,z)\in\Bbb Z^3
\mid n=x^2+y^2+z^2\big\}\big|.$$ Since Legendre and Gauss it is well
known that $r_3(n)>0$ if and only if $n$ is not
 of the form
$4^{\alpha}(8k+7)$, where $\a$ and $k$ are nonnegative integers.
 In 1801 Gauss  (see [D, p.262]) proved that if $n>4$
 is squarefree, then
$$r_3(n)=\cases 24h(-n)&\t{if $n\e 3\mod 8$,}
\\ 12h(-4n)&\t{if $n\e 1,2,5,6\mod 8$,}
\\0&\t{if $n\e 7\mod 8$,}
\endcases\tag 1.11$$
where $h(d)$ is the number of classes consisting of
 primitive binary quadratic forms of discriminant $d$. Let $n\in\Bbb N$ and
$n=2^{\alpha_2} \prod_{i=1}^sp_i^{\alpha_i}$, where $p_1,\ldots,
p_s$ are distinct odd primes. In 1907 Hurwitz (see [D, p.271])
proved that
$$r_3(n^2)=6\prod_{i=1}^s\Big(\f{p_i^{\alpha_i+1}-1}
{p_i-1}-(-1)^{\f{p_i-1}2} \f{p_i^{\alpha_i}-1} {p_i-1}\Big).\tag
1.12$$ For similar formulas concerning $N(a,b,c;n)$ see [B], [CL],
[GPQ], [H], [J], [K1] and [Ye]. In 1862 Liouville (cf. [D, p. 23])
proved that for $a,b,c\in\Bbb N$, $t(a,b,c;n)\ge 1$ for every
$n\in\Bbb N$ if and only if $(a, b, c)=(1, 1, 1), (1, 1, 2), (1, 1,
4), (1, 1, 5), (1, 2, 2), (1, 2, 3)$ or $(1, 2, 4)$. In 1924,
Bell[B] gave transformation relations
 between $N(a,b,c;n)$ and
$r_3(n)$ for
$(a,b,c)=(1,1,2),(1,1,4),(1,1,8),(1,2,2),(1,2,4),(1,2,8),(1,4,4),
(1,4,8),(1,8,8)$. Bell's relations can be easily proved by using
Ramanujan's theta function identities. H$\ddot {\rm u}$rlimann[H]
gave similar results for $(a,b,c)=(1,2,16), (1,8,16)$.

\par Let $a,b,c,d,n\in\Bbb N$. In 2011, the author [S1, Theorem 2.3] found two general relations
between $t(a,b;n)$ and $N(a,b;8n+a+b)$. Recently, using (1.1)-(1.5)
the author and Wang (see [S2],[WS]) revealed
 some new connections between
$t(a,b,c,d;n)$ and $N(a,b,c,d;n)$. They do not need  assuming
$a+b+c+d\le 8$. More recently Yao[Y] confirmed some conjectures
posed by the author in [S2].

\par Let $m,n\in\Bbb N$. In Section 2, using Ramanujan's theta
functions we prove that $$t(1,1,8;n)-\f 13r_3(4n+5) =\cases
2(-1)^{\f{m+1}2}m&\t{if $4n+5=m^2$ for some $m\in\Bbb N$,}
\\0&\t{otherwise}.\endcases\tag 1.13$$
Let  $m\e 1\mod 4$ or $m\e 4\mod 8$. Suppose that
 there is an odd prime
divisor $p$ of $m$ such that
 $\ls{4n+5}p=-1$, where $\sls ap$ is the Legendre
 symbol. Using (1.13) we deduce that
$$t(1,1,8,m;n)=\f 12N(1,1,8,m;8n+10+m),\tag 1.14$$
which confirms  [S2, Conjectures 2.2 and 2.6-2.8].
  We also
show that
$$t(1,3,9;n)=\f 12N(1,3,9;8n+13).\tag 1.15$$

Let $a,b,n\in\Bbb N$ with $2\nmid a$. In Section 3, using
Ramanujan's theta functions we show that

$$t(a,3a,b;n)=\cases 2N(4a,12a,b;8n+4a+b)&\t{if $2\nmid b$,}
\\2N(2a,6a,b/2;4n+2a+b/2)&\t{if $4\mid b-2$,}
\\2N(a,3a, b/4;2n+a+ b/4)-2N(a,3a,b;2n+a+b/4)&\t{if $4\mid b$.}
\endcases\tag 1.16$$

\par Let $a,b,c,d,n\in\Bbb N$.  In Section 5 we
establish many relations between $t(a,b,c,d;n)$ and $N(a,b,c,d;n)$.
For examples, for $2\nmid ab$ we have
$$\align &t(a,a,2b,2b;n)=N(a,a,2b,2b;4n+a+2b),\tag 1.17
\\&t(a,a,b,b;n)=N(a,a,b,b;4n+a+b)-N(a,a,b,b;2n+(a+b)/2),\tag 1.18\endalign$$
$$t(a,3a,c,d;n)=2N(4a,12a,c,d;8n+4a+c+d)\ (\t{if $2\nmid ac$ and $d\e 2,c\pmod 4$}),\tag 1.19$$
 $$\aligned &t(a,3a,2b,c;n)\\&=\f
23\big(N(a,3a,2b,c;8n+4a+2b+c)-N(a,3a,2b,4c;8n+4a+2b+c)\big),\endaligned\tag
1.20$$
$$\aligned &t(a,2a,4a,b;n)
\\&=\f 16\Big(N(a,a,a,2b;16n+14a+2b)
-N\big(a,a,a,2b;4n+3a+\f{a+b}2\big)\Big),
\endaligned\tag 1.21$$
$$t(a,a,2b,4b;n)
=N(a,a,b,2b;4n+a+3b)-N\big(a,a,b,2b;2n+\f{a+3b}2\big),\tag 1.22$$
$$\aligned&t(a,2a,b,2b;n)
\\&=N(a,2a,b,2b;8n+3a+3b)-N(a,2a,b,2b;4n+\f{3(a+b)}2\big)\ (\t{if}\ 4\mid
a-b),\endaligned\tag 1.23$$
$$t(a,3a,9a,d;n)=\f
12\big(N(a,3a,9a,d;8n+13a+d)-N(a,3a,9a,4d;8n+13a+d)\big).\tag 1.24$$
Let $a,b,c,d,m,n\in\Bbb N$ with $2\nmid m$ and $a+b+c\le 7$. We
prove that
$$\aligned t(am,bm,cm,d;n)
&=\f
8{C(a,b,c)}\big(N(am,bm,cm,d;8n+am+bm+cm+d)\\&\q-N(am,bm,cm,4d;8n+am+bm+cm+d)\big).\endaligned\tag
1.25$$
 We also show that
$$\align &t(1,1,1,6;n)
=\f 16\big(N(1,1,1,6;32n+36)-N(1,1,1,6;8n+9)\big),\tag 1.26
\\&t(1,1,1,7;n)=4N(1,1,1,7;4n+5)-2N(1,1,1,7;8n+10),\tag 1.27
\\&t(1,2,6,6;n)=2N(1,2,6,6;8n+15)-N(1,2,6,6;16n+30).\tag 1.28\endalign$$ In
Section 6, based on calculations with Maple we pose some challenging
conjectures.

\section*{2. Formulas for $t(1,1,8;n)$
and $t(1,3,9;n)$}

\par
By (1.1) and (1.2), for $|q|<1$,
$$\aligned\sum_{n=0}^{\infty}r_3(n)q^n&=\varphi(q)^3
=(\phq 4+2q\psq 8)^3
\\&=\phq 4^3+6q\phq 4^2\psq 8+12q^2\phq 4\psq 8^2+8q^3\psq 8^3
\\&=\phq 4^3+6q\phq 4\psq 4^2+12q^2\psq 4^2\psq 8
+8q^3\psq 8^3.
\endaligned\tag 2.1$$
For any $k,s\in\Bbb N$ the power series expansions of
 $\phq {8k}^s$
and $\psq{8k}^s$ are of the form $\sum_{n=0}^{\infty}c_nq^{8n}$.
Thus collecting the terms of the form $q^{4n+2}$ in (2.1) yields
$$\sum_{n=0}^{\infty}r_3(4n+2)q^{4n+2}
=12q^2\psq 4^2\psq 8$$ and so
$$\sum_{n=0}^{\infty}r_3(4n+2)q^n
=12\psi(q)^2\psq 2=12\sum_{n=0}^{\infty}t'(1,1,2;n)q^n.$$ This
yields
$$t(1,1,2;n)=8t'(1,1,2;n)=\f 23r_3(4n+2).\tag 2.2$$
\par By (2.1),
$$\sum_{n=0}^{\infty}
r_3(4n+1)q^{4n+1}=6q\varphi(q^4)\psi(q^4)^2 \ \t{and so} \
\sum_{n=0}^{\infty}r_3(4n+1)q^n=6\varphi(q)\psi(q)^2. \tag 2.3$$
 By (1.5),
 $$\aligned \varphi(q)^2&
=(\phq {16}+2q^4\psq{32}+2q\psq 8)^2
\\&=\phq{16}^2+4q^8\psq{32}^2+4q^4\phq{16}\psq{32}
\\&\qq+4q\phq{16}\psq 8+8q^5\psq 8\psq{32}+4q^2\psq 8^2.
\endaligned\tag 2.4$$
Therefore,
$$\align \sum_{n=0}^{\infty}N(1,1,8;n)q^n
&=\varphi(q)^2\phq 8
=\big(\phq{16}^2+4q^8\psq{32}^2+4q^4\phq{16}\psq{32}
\\&\qq+4q\phq{16}\psq 8+8q^5\psq 8\psq{32}+4q^2
\psq 8^2\big)\phq 8.
\endalign$$
Collecting the terms of the form $q^{8n+2}$ yields
$$\sum_{n=0}^{\infty}N(1,1,8;8n+2)q^{8n+2} =4q^2\psq 8^2\cdot
\phq 8$$ and hence
$$\sum_{n=0}^{\infty}N(1,1,8;8n+2)q^n
=4\varphi(q)\psi(q)^2.\tag 2.5$$ Comparing (2.5) with (2.3) yields
$$N(1,1,8;8n+2)=\f 23r_3(4n+1),\tag 2.6$$
which was first obtained by Bell[B].

\par\q
\pro{Theorem 2.1} Let $n\in\Bbb N$. Then
$$\align &t(1,1,8;n)-\f 12N(1,1,8;8n+10)
\\&=t(1,1,8;n)-\f 13r_3(4n+5)
=\cases 2(-1)^{\f{m+1}2}m&\t{if $4n+5=m^2$ for some $m\in\Bbb N$,}
\\0&\t{otherwise}.\endcases\endalign$$
\endpro
Proof. By (2.6), $N(1,1,8;8n+10)=\f 23r_3(4n+5).$ Set
$s(n)=t(1,1,8;n)-\f 13r_3(4n+5).$ By (2.3) and the fact $r_3(1)=6$,
for $0<|q|<1$ we have
$$\align \sum_{n=0}^{\infty}s(n)q^n
&=\sum_{n=0}^{\infty}t(1,1,8;n)q^n
 -\f 13 \sum_{n=0}^{\infty}r_3(4n+5)q^n
\\& =8\sum_{n=0}^{\infty}t'(1,1,8;n)q^n
- \f 1{3q}
 \sum_{n=0}^{\infty}r_3(4(n+1)+1)q^{n+1}
 \\&=8\psi(q)^2\psq 8-\f 1{3q}\sum_{n=1}^{\infty}
 r_3(4n+1)q^n
 \\&=8\psi(q)^2\psq 8-\f 1{3q}\big(6\psi(q)^2\varphi(q)-
 r_3(1)\big)
 \\&=8\psi(q)^2\psq 8
  -2\f{\psi(q)^2\varphi(q)-1}q
  \\&=\f{2\psi(q)^2(4q\psq 8-\varphi(q))+2}{q}.
  \endalign$$
By [Be2, p.71], $\varphi(-q)=\varphi(q)-4q\psq 8$. Thus,
$$\sum_{n=0}^{\infty}s(n)q^n =2\f{1-\varphi(-q)
\psi(q)^2}q. \tag 2.7$$ It is known that (see [K2, pp.113-114])
$$\psi(q)=\prod_{n=1}^{\infty}\f{(1-q^{2n})^2}{1-q^n}
\qtq{and}
 \varphi(-q)=\prod_{n=1}^{\infty}\f{(1-q^n)^2}
 {1-q^{2n}}.\tag 2.8$$
Thus, appealing to the following Jacobi's identity (see [K2, p.8])
$$\prod_{n=1}^{\infty}(1-q^n)^3=\sum_{n=0}^{\infty}
(-1)^n(2n+1)q^{\f{n(n+1)}2}\q (|q|<1)\tag 2.9$$ we get
$$\align\sum_{n=0}^{\infty}s(n)q^n
&=2\f{1-\varphi(-q)\psi(q)^2}q =\f 2q
\Big(1-\prod_{n=1}^{\infty}\f{(1-q^n)^2}{1-q^{2n}}
\cdot\prod_{n=1}^{\infty}\f{(1-q^{2n})^4}{(1-q^n)^2} \Big)\\&=\f
2q\big(1-\prod_{n=1}^{\infty}(1-q^{2n})^3\big) =\f
2q\Big(1-\sum_{k=0}^{\infty}(-1)^k(2k+1)q^{k(k+1)} \Big)
 \\&=2\sum_{k=1}^{\infty}(-1)^{k+1}(2k+1)
 q^{k^2+k-1}
=2\sum_{k=1}^{\infty}(-1)^{k+1}(2k+1) q^{\f{(2k+1)^2-5}4}.
 \endalign$$
Now comparing the coefficients of $q^n$ on both sides
 yields
 $$s(n)=\cases 2(-1)^{\f{m+1}2}m&\t{if $4n+5=m^2$ for some $m\in\Bbb N$,}
\\0&\t{otherwise}.\endcases$$
This proves the theorem.
\par\q
 \pro{Corollary 2.1} Suppose $n\in\Bbb N$. If
$n\e 0\mod 2$,
 $n\e
0\mod 3$, $n\e 2,3\mod 5$ or $n\e 0,2,3\mod 7$, then
$$t(1,1,8;n)=\f 12N(1,1,8;8n+10).$$
\endpro
Proof. If $2\mid n$, then $4n+5\e 5\mod 8$. If $3\mid n$, then
$4n+5\e 2\mod 3$. If $n\e 2,3\mod 5$, then $4n+5\e 2,3\mod 5$. If
$n\e 0,2,3\mod 7$, then $4n+5\e 3,5,6\mod 7$. Thus, if $n$ satisfies
one of the assumed conditions, then $4n+5$ is not a square and so
$t(1,1,8;n)=\f 12N(1,1,8;8n+10)$ by Theorem 2.1.
\par\q
 \pro{Theorem 2.2}
Let $n\in\Bbb N$. Then $n$ is
 represented by
$\f{x(x-1)}2+\f{y(y-1)}2+8\f{z(z-1)}2$ if and only if
 $4n+5$ is not a
square or $4n+5$ has a prime
 divisor of the form $4k+3$.
\endpro
Proof. Since $4n+5\not\e 7\mod 8$ we have $r_3(4n+5)>0$. Thus, when
$4n+5$ is not a square we have $t(1,1,8;n)= \f 13 r_3(4n+5)>0$ by
Theorem 2.1. Now assume that $4n+5=(2a+1)^2$ for
 $a \in\Bbb N$. By Theorem 2.1,
$$t(1,1,8;n)= \f 13 r_3(4n+5)+(-1)^{a+1}(4a+2).$$
If $a$ is odd, then $2a+1\e 3\mod 4$,
 $t(1,1,8;n) =\f 13
r_3(4n+5) +4a+2>0$ and $4n+5$ has a prime
 divisor of the form $4k+3$.
 Now assume that $a=2m$ for $m\in\Bbb N$.
Then $4n+5=(4m+1)^2$.
  By Theorem 2.1,
$$t(1,1,8;n)= \f 13 r_3(4n+5)-(8m+2).$$
By (1.12), if $4m+1=p_1^{\a_1}\cdots p_s^{\a_s}$, where
$p_1,\ldots,p_s$ are distinct primes, then
$$\align r_3(4n+5)&=r_3((4m+1)^2)=6\prod_{i=1}^s
\Big(\f{p_i^{\a_i+1}-1}{p_i-1}-(-1)^{\f{p_i-1}2}
\f{p_i^{\a_i}-1}{p_i-1}\Big)
\\&\q\ge 6\prod_{i=1}^s
\Big(\f{p_i^{\a_i+1}-1}{p_i-1} -\f{p_i^{\a_i}-1}{p_i-1}\Big)
=6p_1^{\a_1}\cdots p_s^{\a_s}=6(4m+1).
\endalign$$
Moreover, the strict inequality holds if and only
 if $4m+1$ has a
prime divisor of the form $4k+3$. Hence, $t(1,1,8;n)=
 \f 13r_3(4n+5)-2(4m+1)>0$ if and only if
 $4m+1$ has a prime
divisor of the form $4k+3$. This completes the proof.

 \par\q
\pro{Theorem 2.3} Let $m,n\in\Bbb N$
 with $m\e 1\mod 4$ or
$m\e 4\mod 8$. Suppose that
 there is an odd prime
divisor $p$ of $m$ such that
 $\ls{4n+5}p=-1$. Then
$$t(1,1,8,m;n)=\f 12N(1,1,8,m;8n+10+m).$$
\endpro
Proof. Suppose that $p$ is an odd prime divisor of $m$ with
 $\ls{4n+5}p=-1$. For $w\in \Bbb Z$ we see that
$$\Ls{4(n-m\f{w(w-1)}2)+5}p
=\Ls{4n+5}p=-1.$$ Hence $4(n-m\f{w(w-1)}2)+5$ is not a square. Now,
from Theorem 2.1 we derive that
$$\align t(1,1,8,m;n)
&=\sum_{w\in\Bbb Z}t \big(1,1,8;
 n-m\f{w(w-1)}2\big)
\\&=\f 12\sum_{w\in\Bbb Z}N
(1,1,8;8n+10-m\cdot 4w(w-1))
\\&=\f 12\sum_{w\in\Bbb Z}N
(1,1,8;8n+10+m-m(2w-1)^2).
\endalign$$
 Since $a^2\e 0,1 \mod
4$ and $a^2\e 0,1,4\mod 8$ for any $a\in\Bbb Z$, we see that
$x^2+y^2
 \not\e 3\mod 4$ and $x^2+y^2\not\e 6\mod 8$
 for any $x,y\in\Bbb Z$.
  If $m\e 1\mod 4$ and
  $8n+10+m-m(2w)^2=x^2+y^2+8z^2$ for some
$x,y,z,w\in\Bbb Z$, then $x^2+y^2\e 10+m\e 3\mod 4$.
 This is
impossible. If $m\e 4\mod 8$ and
  $8n+10+m-m(2w)^2=x^2+y^2+8z^2$ for some
$x,y,z,w\in\Bbb Z$, then $x^2+y^2\e 10+m\e 6\mod 8$.
 This is also
impossible.
 Hence, for $m\e 1\mod 4$ or $m\e 4\mod 8$,
$$\align t(1,1,8,m;n)
&=\f 12\sum_{w\in\Bbb Z}N (1,1,8;8n+10+m -m(2w-1)^2)
\\&=\f
12\sum_{w\in\Bbb Z}N (1,1,8;8n+10+m-mw^2)
\\&=\f 12 N(1,1,8,m;8n+10+m).
\endalign$$
This proves the theorem.
\par\q \pro{Corollary 2.2 ([S2,
Conjectures 2.6 and 2.8])} Let $n\in\Bbb N$. Then
$$t(1,1,5,8;n)=\f 12N(1,1,5,8;8n+15)\qtq{for}n\e 2,3
\mod 5$$ and
$$t(1,1,8,13;n)=\f 12N(1,1,8,13;8n+23)
\qtq{for}n\e 0,4,7,8,9,10 \mod {13}.$$
\endpro
Proof. Putting $m=5,13$ in Theorem 2.3 yields the result.
\par\q
\par{\bf Remark 2.1} By Theorem 2.3, for $n\e 0\mod 3$ we have
$\sls{4n+5}3=-1$ and so $t(1,1,8,9;n)=\f 12N(1,1,8,9;8n+19)$ and
 $t(1,1,8,12;n)=\f 12N(1,1,8,12;8n+22)$,
which were conjectured by the author in [S2,
 Conjectures 2.2 and
2.7] and first confirmed by Yao in [Y].
\par\q

\par\q
\par For $a,b,c,n\in\Bbb N$ it is clear that
$$\align &n=a\f{x(x-1)}2+b\f{y(y-1)}2+c\f{z(z-1)}2
\\&\iff 8n+a+b+c=a(2x-1)^2+b(2y-1)^2+c(2z-1)^2.
\endalign$$
Thus,
$$t(a,b,c;n)=\big|\big\{(x,y,z)\in\Bbb Z^3\mid 8n+a+b+c
=ax^2+by^2+cz^2,\ 2\nmid xyz\big\}\big|.\tag 2.10$$
\par\q
\pro{Lemma 2.1 ([S2, Lemma 2.3])} For $|q|<1$ we have
$$\align\varphi(q)\phq 3&=
\phq {16}\phq{48}+4q^{16}\psq{32}\psq{96}+2q \phq {48}\psq
8+2q^3\phq{16}\psq{24} \\&\q+6q^4\psq 8\psq{24}+4q^{13}\psq
8\psq{96}+4q^7\psq{24}\psq{32}.\endalign$$
\endpro
 \pro{Theorem 2.4}  For $n\in\Bbb N$ we have
$$t(1,3,9;n)=t(1,3,27;3n+1)=\f 12N(1,3,9;8n+13).$$
\endpro
 Proof.
 By (1.5), (1.7) and Lemma 2.1,
$$\aligned &\sum_{n=0}^{\infty}N(1,3,9;n)q^n
=\varphi(q)\phq{3}\phq{9}
\\&=\big(\phq {16}\phq{48}+4q^{16}\psq{32}\psq{96}
+2q
 \phq {48}\psq
{8}+2q^{3}\phq{16}\psq{24}
\\&\q+6q^{4}\psq {8}\psq{24}+4q^{13}\psq
{8}\psq{96}+4q^{7}\psq{24}\psq{32}\big)
\\&\q\times\big(\phq{144}+2q^{36}\psq{288}+2q^9
\psq{72}\big).
\endaligned\tag 2.11$$
For any $k,s\in\Bbb N$ the power series expansions of
 $\phq{8k}^s$
and $\psq{8k}^s$ are of the form $\sum_{n=0}^{\infty}c_nq^{8n}$.
Thus, collecting the terms of the form
 $q^{8n+13}$ in (2.11) and then applying (1.4) and
 (1.8) we
deduce that
$$\align &\sum_{n=0}^{\infty}N(1,3,9;8n+13)q^{8n+13}
\\&=4q^{13}\psq
{8}\psq{96}\cdot \phq{144} +2q
 \phq {48}\psq
{8}\cdot 2q^{36}\psq{288} +6q^{4}\psq {8}\psq{24}\cdot 2q^9 \psq{72}
\\&=4q^{13}\psq 8\big(\phq{144}\psq{96}+q^{24}
\phq{48}\psq{288}\big)+12q^{13}\psq 8\psq {24}\psq{72}
\\&=16q^{13}\psq{8}\psq{24}\psq{72}
\\&=16q^{13}\sum_{n=0}^{\infty}t'(1,3,9;n)q^{8n}.
\endalign$$
 Now comparing the coefficients of $q^{8n+13}$ on
 both sides yields
$$N(1,3,9;8n+13)=16t'(1,3,9;n)=2t(1,3,9;n).$$
If $8(3n+1)+31=24n+39=x^2+3y^2+27z^2$ for some odd integers $x,y$
and $z$, then clearly $3\mid x$ and so $8n+13=\f 13((3x)^2+
3y^2+27z^2)=3x^2+y^2+9z^2$ for odd integers $x,y$ and $z$. Thus,
applying (2.10) we have
$$\align t(1,3,27;3n+1)&=\big|\big\{(x,y,z)
\in\Bbb Z^3\bigm| 8(3n+1)+31=x^2+3y^2+27z^2, \ 2\nmid xyz\big\}\big|
\\&=\big|\big\{(x,y,z)\in\Bbb Z^3\bigm|
8n+13=3x^2+y^2+9z^2,\ 2\nmid xyz\big\}\big|
\\&=t(1,3,9;n).\endalign$$
This completes the proof.
\par\q
\par{\bf Remark 2.2} One can similarly prove that
$$t(1,1,3;n)=\f 12N(1,1,3;8n+5)\qtq{and}
t(1,3,3;n)=\f 12N(1,3,3;8n+7),$$ which can be deduced from (1.9).
\par\q
\section*{3. Formulas for $t(a,3a,b;n)$}
\par\q
\par By (1.5), (1.7) and Lemma 2.1, for
$a,b\in\Bbb N$ with $2\nmid a$ we have
$$\aligned &\sum_{n=0}^{\infty}N(a,3a,2b;n)q^n
=\phq{a}\phq{3a}\phq{2b}
\\&=\big(\phq {16a}\phq{48a}+4q^{16a}\psq{32a}\psq{96a}+2q^a
 \phq {48a}\psq
{8a}\\&\q+2q^{3a}\phq{16a}\psq{24a} +6q^{4a}\psq
{8a}\psq{24a}+4q^{13a}\psq
{8a}\psq{96a}\\&\q+4q^{7a}\psq{24a}\psq{32a}\big)
\big(\phq{8b}+2q^{2b}\psq{16b}\big).
\endaligned\tag 3.1$$
\pro{Theorem 3.1} Let $a,b\in\{1,3,5,\ldots\}$.
 For $n\in\Bbb N$ we
have
$$t(a,3a,2b;n)=\f 23N(a,3a,2b;8n+4a+2b).$$
\endpro
Proof. For any $k,s\in\Bbb N$ the power series expansions of
$\phq{8k}^s$ and $\psq{8k}^s$ are of the form
$\sum_{n=0}^{\infty}c_nq^{8n}$.
 Since $4a+2b\e 2\mod 4$, collecting the terms of the form $q^{8n+4a+2b}$ in (3.1)
yields
$$\sum_{n=0}^{\infty}N(a,3a,2b;8n+4a+2b)q^{8n+4a+2b}
=6q^{4a}\psq{8a}\psq{24a}\cdot 2q^{2b}\psq{16b}.$$
  Replacing $q$ with $q^{1/8}$ gives
$$\align &\sum_{n=0}^{\infty}N(a,3a,2b;8n+4a+2b)q^n
\\&=12\psq{a}\psq{3a}\psq{2b}
=12\sum_{n=0}^{\infty}t'(a,3a,2b;n)q^n
=\f{12}8\sum_{n=0}^{\infty}t(a,3a,2b;n)q^n.
\endalign$$
Now comparing the coefficients of $q^n$ on both sides yields the
result.
\par\q

\pro{Theorem 3.2} Let $a\in\{1,3,5,\ldots\}$ and $m\in\Bbb N$. For
$n\in\Bbb N$ we have
$$t(a,3a,8m;n)=\f 23N(a,3a,8m;8n+4a+8m)
-2N(a,3a,8m;2n+a+2m).$$
\endpro
Proof. Set $b=4m$. Collecting the terms of the form $q^{8n+4a}$ in
(3.1) we deduce that
$$\sum_{n=0}^{\infty}N(a,3a,2b;8n+4a)q^{8n+4a}
=6q^{4a}\psq{8a}\psq{24a} (\phq{8b}+2q^{2b}\psq{16b}).$$ Replacing
$q$ with $q^{1/8}$ we obtain
$$\sum_{n=0}^{\infty}N(a,3a,8m;8n+4a)q^n
=6\psq{a}\psq{3a} (\phq{4m}+2q^m\psq{8m}). \tag 3.2$$ On the other
hand, using (1.5) and (1.7) we see that
$$\align &\sum_{n=0}^{\infty}N(a,3a,8m;n)q^n
\\&=\phq a\phq{3a}\phq{8m}
=(\phq{4a}+2q^a\psq{8a})(\phq{12a}+2q^{3a}\psq{24a})\phq{8m}.
\endalign$$
Collecting the terms of the form $q^{2n+a}$ and then applying (1.4)
we get
$$\align &\sum_{n=0}^{\infty}N(a,3a,8m;2n+a)q^{2n+a}
\\&=(2q^a\psq{8a}\phq{12a}+2q^{3a}\psq{24a}\phq{4a})\phq{8m}
\\&=2q^a\psq{2a}\psq{6a}\phq{8m}\endalign$$
and so
$$\sum_{n=0}^{\infty}N(a,3a,8m;2n+a)q^n
=2\psq{a}\psq{3a}\phq{4m}.$$ This together with
 (3.2) yields
$$\align &\sum_{n=0}^{\infty}(N(a,3a,8m;8n+4a)-3N(a,3a,8m;2n+a))q^n
\\&=12q^m\psq a\psq{3a}\psq{8m}=12q^m\sum_{n=0}^{\infty}
t'(a,3a,8m;n)q^n
\\&=\f{12}8q^m\sum_{n=0}^{\infty} t(a,3a,8m;n)q^n.
\endalign$$
Now comparing the coefficients of $q^{m+n}$ on both sides gives
 the result.

\par\q
 \pro{Theorem 3.3} Let
$a\in\{1,3,5,\ldots\}$ and $m\in\{0,1,2,\ldots\}$. For $n\in\Bbb N$
we have
$$\align &t(a,3a,8m+4;n)
\\&=\cases \f 23N(a,3a,8m+4;8n+4a+8m+4)
\q \t{if $n\e \f{a-1}2+m\mod 2$,} \\\f
23\big(N(a,3a,8m+4;8n+4a+8m+4) -N(a,3a,8m+4;2n+a+2m+1)\big)
\\\qq\qq\qq\qq\qq\qq\qq\qq\ \;\t{if $n\not\e \f{a-1}2+m\mod 2$}.
\endcases\endalign$$
\endpro
Proof. Set $b=4m+2$. Collecting the terms
 of the form $q^{8n}$ in
(3.1) we deduce that
$$\align &\sum_{n=0}^{\infty}N(a,3a,2b;8n)q^{8n}
\\&=\big(\phq {16a}\phq{48a}+4q^{16a}\psq{32a}\psq{96a}\big)\phq{8b}
+6q^{4a}\psq{8a}\psq{24a}\cdot 2q^{2b}\psq{16b} .\endalign$$
Replacing $q$ with $q^{1/8}$ we obtain
$$\aligned \sum_{n=0}^{\infty}N(a,3a,8m+4;8n)q^n
&=\big(\phq {2a}\phq{6a}+4q^{2a}\psq{4a}\psq{12a}\big)\phq{4m+2}
\\&\q+12q^{m+(a+1)/2}\psq{a}\psq{3a}\psq{8m+4}.
\endaligned\tag 3.3$$
 On the other
hand, using (1.5) we see that
$$\align &\sum_{n=0}^{\infty}N(a,3a,8m+4;n)q^n
\\&=\phq a\phq{3a}\phq{8m+4}
=(\phq{4a}+2q^a\psq{8a})(\phq{12a}+2q^{3a}\psq{24a})\phq{8m+4}.
\endalign$$
Collecting the even powers of $q$ we get
$$\sum_{n=0}^{\infty}N(a,3a,8m+4;2n)q^{2n}
=(\phq{4a}\phq{12a}+4q^{4a}\psq{8a}\psq{24a})\phq{8m+4}$$
 and so
$$\sum_{n=0}^{\infty}N(a,3a,8m+4;2n)q^n
=(\phq{2a}\phq{6a}+4q^{2a}\psq{4a}\psq{12a})\phq{4m+2}.$$
 This together with (3.3) yields
$$\align &\sum_{n=0}^{\infty}(N(a,3a,8m+4;8n)-N(a,3a,8m+4;2n))q^n
\\&=12q^{m+(a+1)/2}\psq{a}\psq{3a}\psq{8m+4}
=12q^{m+(a+1)/2}
 \sum_{n=0}^{\infty}t'(a,3a,8m;n)q^n
\\&=\f 32q^{m+(a+1)/2}\sum_{n=0}^{\infty} t(a,3a,8m+4;n)q^n.\endalign$$
Comparing the coefficients of $q^{m+(a+1)/2+n}$ on both sides
 yields
$$\align
&t(a,3a,8m+4;n)
\\&=\f 23\big(N(a,3a,8m+4;8n+4a+8m+4)
-N(a,3a,8m+4;2n+a+2m+1)\big).\endalign$$ Now assume $n\e
\f{a-1}2+m\mod 2$. Then $2n+2m+a+1\e
 a-1+2m+2m+a+1\e 2a\mod 4$. If
$2n+2m+a+1=ax^2+3ay^2+(8m+4)z^2$ for some $x,y,z\in\Bbb Z$, we must
have $a(x^2+3y^2)\e 2n+2m+a+1\e 2a\mod 4$ and so $x^2+3y^2 \e 2\mod
4$. If $2\mid x-y$, then $4\mid x^2+3y^2$. If $2\nmid x-y$, then
$x^2+3y^2$ is odd. Thus, $x^2+3y^2\not\e 2\mod 4$ and we get a
contradiction. Therefore $N(a,3a,8m+4;2n+a+2m+1)=0$. This completes
the proof.

\pro{Theorem 3.4} Let $a,b,n\in\Bbb N$ with $2\nmid a$. Then
$$t(a,3a,b;n)=\cases 2N(4a,12a,b;8n+4a+b)&\t{if $2\nmid b$,}
\\2N(2a,6a,b/2;4n+2a+b/2)&\t{if $4\mid b-2$,}
\\2N(a,3a, b/4;2n+a+ b/4)-2N(a,3a,b;2n+a+b/4)&\t{if $4\mid b$.}
\endcases$$
\endpro
Proof. Clearly the result is equivalent to
$$t(a,3a,b;n)=2(N(4a,12a,b;8n+4a+b)-N(4a,12a,4b;8n+4a+b)).$$
By (1.2) and (1.4), $\phq b-\phq{4b}=2q^b\psq{8b}$ and
$$\align&
\phq{4a}\phq{12a}\\&=(\phq{16a}+2q^{4a}\psq{32a})(\phq{48a}+2q^{12a}\psq{96a})
\\&=\phq{16a}\phq{48a}+4q^{16a}\psq{32a}\psq{96a}+2q^{4a}\phq{48a}\psq{32a}
+2q^{12a}\phq{16a}\psq{96a}
\\&=\phq{16a}\phq{48a}+4q^{16a}\psq{32a}\psq{96a}+2q^{4a}\psq{8a}\psq{24a}.
\endalign$$
Hence
$$\align &\sum_{n=0}^{\infty}(N(4a,12a,b;n)-N(4a,12a,4b;n))q^n
\\&=\phq{4a}\phq{12a}(\phq b-\phq{4b})
\\&=2q^b\psq{8b}(\phq{16a}\phq{48a}+4q^{16a}\psq{32a}\psq{96a}+2q^{4a}\psq{8a}\psq{24a}).
\endalign$$
Therefore
$$\sum\Sb n=0\\n\e 4a+b\mod 8\endSb^{\infty}(N(4a,12a,b;n)-N(4a,12a,4b;n))q^n
=4q^{4a+b}\psq{8a}\psq{24a}\psq{8b}$$ and so
$$\sum_{n=0}^{\infty}(N(4a,12a,b;8n+4a+b)-N(4a,12a,4b;8n+4a+b))q^{8n}=4\psq{8a}\psq{24a}\psq{8b}.$$
This yields
$$\align &\sum_{n=0}^{\infty}(N(4a,12a,b;8n+4a+b)-N(4a,12a,4b;8n+4a+b))q^n
\\&=4\psq a\psq {3a}\psq b=\f 12\sum_{n=0}^{\infty}t(a,3a,b;n)q^n.\endalign$$
Comparing the coefficients of $q^n$ on both sides gives the result.
\par
\par Comparing Theorems 3.1-3.3 with Theorem 3.4 we deduce the
following result.
 \pro{Corollary 3.1} Suppose $a,b,n\in\Bbb N$ with $2\nmid a$. Then
$$\align &N(a,3a,2b;8n+4a+2b)\\&=\cases 3N(2a,6a,b;4n+2a+b)&\t{if $2\nmid b$,}
\\3N(a,3a,b/2;2n+a+b/2)&\t{if $4\mid b$,}
\\3N(a,3a,b/2;2n+a+b/2)-2N(a,3a,2b;2n+a+b/2)&\t{if $4\mid b-2$.}
\endcases\endalign$$
\endpro

\section*{4. General relations between $t(a,b,c;n)$ and $N(a,b,c;
8n+a+b+c)$}
\par\q
\pro{Theorem 4.1} Let $a,b,c\in\Bbb N$ with $2\nmid ab$, $4\mid a-b$
and $4\mid c-2$. For $n\in\Bbb N$ we have
$$t(a,b,c;n)=N(a,b,c;8n+a+b+c)-N(a,b,c;2n+(a+b+c)/4).$$
\endpro
Proof. Suppose $8n+a+b+c=ax^2+by^2+cz^2$ for some $x,y,z\in\Bbb Z$.
Since  $a\e b\e \pm 1\mod 4$, $c\e 2\mod 4$ and $a+b+c\e 0\mod 4$,
we see that either $x\e y\e z\e 0\mod 2$ or $x\e y\e z\e 1\mod 2$.
Hence appealing to (2.10) we deduce that
$$\align &t(a,b,c;n)
\\&=\big|\big\{(x,y,z)\in\Bbb Z^3\bigm| 8n+a+b+c=ax^2
+by^2+cz^2,\ 2\nmid xyz\big\}\big|
\\&=N(a,b,c;8n+a+b+c)
\\&\q-\big|\big\{(x,y,z)\in\Bbb Z^3\bigm| 8n+a+b+c=ax^2 +by^2+cz^2,\ x\e
y\e z\e 0\mod 2\big\}\big|
\\&=N(a,b,c;8n+a+b+c)
\\&\q-\big|\big\{(x,y,z)\in\Bbb Z^3\bigm| 8n+a+b+c=a(2x)^2
+b(2y)^2+c(2z)^2\big\}\big|
\\&=N(a,b,c;8n+a+b+c)-N(a,b,c;2n+(a+b+c)/4).
\endalign$$
This proves the theorem.
\par\q

\pro{Theorem 4.2} Suppose that $a,b,c \in\Bbb N$ with $2\nmid ab$
and $4\mid a-b$. If $c\e a\mod 4$ or $c\e 4\mod 8$, then
$$t(a,b,c;n)=N(a,b,c;8n+a+b+c). $$
\endpro
Proof.
 Assume $c\e a\mod 4$ and
  $8n+a+b+c=ax^2+by^2+cz^2$ for some
  $x,y,z\in\Bbb Z$. If $2\mid z$, then $3a\e
   8n+a+b+c\e ax^2+by^2\e a(x^2+y^2)\mod 4$.
    Since $x^2,y^2
   \e 0,1\mod 4$, we must have
    $x^2+y^2\not\e 3\mod 4$ and get a contradiction.
     Hence $2\nmid z$. Then
     $a(x^2+y^2)\e
     ax^2+by^2=8n+a+b+c-cz^2\e a+b\e 2a\mod 4$. That is,
     $x^2+y^2\e 2\mod 4$. This implies that
     $2\nmid xy$.
\par Now assume
 $c\e 4\mod 8$ and $8n+a+b+c=ax^2+by^2+cz^2$ for some
  $x,y,z\in\Bbb Z$. Then $a(x^2+y^2)\e ax^2+by^2
  =8n+a+b+c-cz^2\e a+b\e 2\mod 4$. This implies $2\nmid xy$ and so
  $cz^2=8n+a+b+c-ax^2-by^2\e a+b+c-a-b=c\mod 8$.
  Hence $z$ is also odd.
\par By the above and (2.10),
for $c\e a\mod 4$ or $c\e 4\mod 8$,
$$\align t(a,b,c;n)
&=\big|\big\{(x,y,z)\in\Bbb Z^3\bigm| 8n+a+b+c=ax^2 +by^2+cz^2,\
2\nmid xyz\big\}\big|
\\&=N(a,b,c;8n+a+b+c).\endalign$$
This proves the theorem.
\par\q
\pro{Theorem 4.3} Let $a,b,c\in\Bbb N$ with $2\nmid
 a$, $2\mid b$, $2\mid c$, $8\nmid b$, $8\nmid c$
  and $8\nmid b+c$. Then $t(a,b,c;n)=N(a,b,c;8n+a+b+c)$.
\endpro
Proof. Suppose $8n+a+b+c=ax^2+by^2+cz^2$ for
 some $x,y,z\in\Bbb Z$.
Then clearly $2\nmid x$ and so $by^2+cz^2 =8n+a+b+c-ax^2\e b+c\mod
8$. If $2\mid y$, since $8\nmid b+c$ we have
 $2\nmid z$ and so $c\e cz^2\e
by^2+cz^2\e b+c\mod 8$. This contradicts the assumption $8\nmid b$.
Hence $2\nmid y$. Similarly, $2\nmid z$. Now applying (2.10) yields
the result.
\par\q
 \pro{Lemma 4.1} Let $a_1,\ldots,a_k,d,m,n\in\Bbb N$ with $k\ge 2$ and $2\nmid m$.
If
$$t(a_1,\ldots,a_k;s)=cN(a_1,\ldots,a_k;8s+a_1+\cdots+a_k)\qtq{for}
s=0,1,2,\ldots,$$ then
$$\align t(a_1m,\ldots, a_km,d;n)
&=c\big(N(a_1m,\ldots,a_km,d;8n+(a_1+\cdots+a_k)m+d)
\\&\q-N\big(a_1m,\ldots,a_km,4d;8n+(a_1+\cdots+a_k)m+d)\big).\endalign$$
\endpro
 Proof. Clearly
$$\align &t(a_1m,\ldots, a_km,d;n)
\\&=\sum_{w\in\Bbb Z}t(a_1m,\ldots, a_km;n-dw(w-1)/2)
\\&=\sum_{w\in\Bbb Z,m\mid n-d\f{w(w-1)}2}
t\big(a_1,\ldots,a_k;\f{n-dw(w-1)/2}m\big)
\\&=c\sum_{w\in\Bbb Z,m\mid
n-d\f{w(w-1)}2}N\big(a_1,\ldots,a_k;8\f{n-dw(w-1)/2)}m
+a_1+\cdots+a_k\big)
\\&=c\sum_{w\in\Bbb Z}N(a_1m,\ldots, a_km;
8n+a_1m+\cdots+a_km+d-d(2w-1)^2)
\\&=c\Big(\sum_{w\in\Bbb Z}N(a_1m,\ldots, a_km;
8n+a_1m+\cdots+a_km+d-dw^2)\\&\qq-\sum_{w\in\Bbb Z}N(a_1m,\ldots,
a_km;8n+a_1m+\cdots+a_km+d-d(2w)^2)\Big)
\\&=c\big(N(a_1m,\ldots, a_km,d;8n+a_1m+\cdots+a_km+d)
\\&\qq-N(a_1m,\ldots, a_km,4d;8n+a_1m+\cdots+a_km+d)\big).
\endalign$$
This proves the lemma.
\par\q
\par Now we present the following general theorem.

\pro{Theorem 4.4} Let $a_1,\ldots,a_k,d,m,n\in\Bbb N$ with $2\nmid
m$, $k\ge 2$ and $a_1+\cdots+a_k\le 7$. Then
$$\align t(a_1m,\ldots, a_km,d;n)
&=\f {2^k}C\big(N(a_1m,\ldots,a_km,d;8n+(a_1+\cdots+a_k)m+d)
\\&\q-N(a_1m,\ldots,a_km,4d;8n+(a_1+\cdots+a_k)m+d)\big),\endalign$$
where
$$C=2^k+2^{k-1}\b{i_1}4+2^{k-2}
i_1(i_1-1)i_2+2^{k-1}i_1i_3$$ and $i_j$ is the number of elements in
$\{a_1,\ldots,a_k\}$ which are equal to $j$.
\endpro
Proof. By (1.6) and (1.9),
$$t(a_1,\ldots,a_k;n)=\f {2^k}CN(a_1,\ldots,a_k;8n+a_1
+\cdots+a_k).$$ Hence applying Lemma 4.1 yields the result.
\par\q
\pro{Corollary 4.1} Let $a,b,n\in\Bbb N$ with $2\nmid a$. Then
$$t(a,3a,b;n)=\f 23\big(N(a,3a,b;8n+4a+b)-N(a,3a,4b;8n+4a+b)\big).$$
\endpro

\pro{Lemma 4.2 ([S1, Theorem 2.3])} Suppose $a,b,n\in\Bbb N$,
$8\nmid a$, $8\nmid b$ and $4\nmid a+b$. Then
$t(a,b;n)=N(a,b;8n+a+b)$.\endpro
\par\q
 \pro{Theorem 4.5} Suppose
$a,b,c,n\in\Bbb N$, $8\nmid a$, $8\nmid b$ and $4\nmid a+b$. Then
$$t(a,b,c;n)=N(a,b,c;8n+a+b+c)-N(a,b,4c;8n+a+b+c).$$
\endpro
Proof. By Lemma 4.2, $t(a,b;n)=N(a,b;8n+a+b)$. Now applying Lemma
4.1 (with $m=1$) gives the result.
\par\q
\pro{Lemma 4.3 ([S1, Theorem 2.3])} Suppose $a,b,n\in\Bbb N$,
$2\nmid a$, $8\mid b-4$ and $4\mid a+\f b4$. Then
$t(a,b;n)=N(a,b/4;8n+a+b)$.\endpro
\par\q
 \pro{Theorem 4.6} Suppose
$a,b,c,n\in\Bbb N$, $2\nmid a$, $8\mid b-4$ and $4\mid a+\f b4$.
Then
$$t(a,b,c;n)=N\big(a,\f b4,c;8n+a+b+c\big)-N\big(a,\f b4,4c;8n+a+b+c\big).$$
\endpro
Proof. Using Lemma 4.3 we see that
$$\align &t(a,b,c)
\\&=\sum_{z\in\Bbb Z}t(a,b;n-cz(z-1)/2)
\\&=\sum_{z\in\Bbb Z}N\big(a,\f b4;8(n-cz(z-1)/2)+a+b\big)
\\&=\sum_{z\in\Bbb Z}N\Big(a,\f
b4;8\big(n-c\f{z(z-1)}2\big)+a+b\Big)
\\&=\sum_{z\in\Bbb Z}N\Big(a,\f
b4;8n+a+b+c-c(2w-1)^2\Big)
\\&=\sum_{z\in\Bbb Z}N\Big(a,\f
b4;8n+a+b+c-cw^2\Big)-\sum_{z\in\Bbb Z}N\Big(a,\f
b4;8n+a+b+c-c(2w)^2\Big)
\\&=N\big(a,\f b4,c;8n+a+b+c\big)-N\big(a,\f b4,4c;8n+a+b+c\big).
\endalign$$
This proves the theorem.

\par\q
\pro{Theorem 4.7} Suppose $n\in\Bbb N$. Then $n$ is represented by
$\f{x(x-1)}2+\f{y(y-1)}2 +9\f{z(z-1)}2$ if and only if $n\not\e
5,8\mod 9$.
\endpro
Proof. By Theorem 4.2, $t(1,1,9;n)=N(1,1,9;8n+11)$. It is well known
that $8n+11=x^2+y^2+z^2$ for some $x,y,z\in\Bbb Z$.
 Assume $n\e 0,1\mod
3$. Then $3\nmid 8n+11$. Since $m^2\e 0,1\mod 3$ for $m\in\Bbb Z$ we
must have $3\mid xyz$ and so $8n+11$ is represented by
$x^2+y^2+9z^2$. Hence $n$ is represented by $\f{x(x-1)}2+\f{y(y-1)}2
+9\f{z(z-1)}2$. For $n\e 2\mod 9$ we have
 $9\mid 8n+11$ and
$\f{8n+11}9\e 3\mod 8$. Hence $\f{8n+11}9=x^2+y^2+z^2$ for some
$x,y,z\in\Bbb Z$. This yields $8n+11=(3x)^2+(3y)^2+9z^2$. Therefore
$t(1,1,9;n)=N(1,1,9;8n+11)>0$. Finally we assume $n\e 5,8\mod 9$. If
  $8n+11=x^2+y^2+9z^2$ for some integers
$x,y$ and $z$, since $m^2\e 0,1,4,7\mod 9$ for $m\in\Bbb Z$
 we see
that $x^2+y^2\e 0,1,2,4,5,7,8\mod 9$ and so
$8n+11=x^2+y^2+9z^2\not\e 3,6\mod 9$. This yields $n\not\e 5,8\mod
9$, which contradicts the assumption. Hence
$t(1,1,9;n)=N(1,1,9;8n+11)=0$ for $n\e 5,8\mod 9$.
 Putting the above together
proves the theorem.

\section*{5. Some relations between $t(a,b,c,d;n)$ and $N(a,b,c,d;n)$}

\pro{Theorem 5.1} Suppose $a,b,c,d,n\in\Bbb N, 2\nmid abc$ and $a\e
b\e c\pmod 4.$ Then
$$\align &t(a,b,c,d;n)
\\&=N(a,b,c,d;8n+a+b+c+d)-N(a,b,c,4d;8n+a+ b+c+d).\endalign$$
\endpro
Proof. By Theorem 4.2, $t(a,b,c;k)=N(a,b,c;8k+a+b+c)$ for any
$k\in\Bbb N.$ Now applying Lemma 4.1  (with $m=1$) yields the
result.
\pro{Corollary 5.1} Let $a,b,c,d,n\in\Bbb N$ with $a\e b \e
c\e \pm 1\mod 4$ and $d\e 4\mod 8$. Then
$$t(a,b,c,d;n)=N(a,b,c,d;8n+a+b+c+d).$$
\endpro
Proof. By Theorem 5.1,
$$t(a,b,c,d;n)=N(a,b,c,d;8n+a+b+c+d)-N(a,b,c,4d;8n+a+b+c+d).$$
If $8n+a+b+c+d=ax^2+by^2+cz^2+4dw^2$ for $x,y,z,w\in\Bbb Z,$ we see
that $$a(x^2+y^2+z^2)\e ax^2+by^2+cz^2\e a+b+c\e 3a\pmod 4$$ and so
$x^2+y^2+z^2\e 3\pmod 4.$ This yields $2\nmid xyz$ and so
$$ax^2+by^2+cz^2\e a+b+c\not\e 8n+a+b+c+d-4dw^2\pmod 8,$$ which is a
contradiction. Hence, $N(a,b,c,4d;8n+a+b+c+d)=0$ and the result
follows.

\pro{Theorem 5.2} Suppose $a,b,c,d,n\in\Bbb N$, $2\nmid abcd$ and
$a\e b\e c\e d\pmod 4.$ Then
$$t(a,b,c,d;n)=N\big(a,b,c,d;8n+a+b+c+d\big)-
N\Big(a,b,c,d;2n+\f{a+b+c+d}2\Big).$$
\endpro
Proof. By Theorem 5.1,
$$t(a,b,c,d;n)=N(a,b,c,d;8n+a+b+c+d)-N(a,b,c,4d;8n+a+b+c+d).$$
If $8n+a+b+c+d=ax^2+by^2+cz^2+4dw^2$ for some $x,y,z,w\in\Bbb Z,$
then $$a(x^2+y^2+z^2)\e ax^2+by^2+cz^2=8n+a+b+c+d-4dw^2\e0\pmod 4$$
and so $4\mid x^2+y^2+z^2.$ This implies that $2\mid x, 2\mid y$ and
$2\mid z.$ Hence,
$$N(a,b,c,4d;8n+a+b+c+d)=N\Big(a,b,c,d;2n+\f{a+b+c+d}4\Big).$$
So the result follows.

\pro{Theorem 5.3} Suppose $a,b,c,d,n\in\Bbb N, 2\nmid a,2\mid b,
2\mid c, 8\nmid b, 8\nmid c$ and $8\nmid b+c.$ Then
$$\align &t(a,b,c,d;n)\\&=N(a,b,c,d;8n+a+b+c+d)-N(a,b,c,4d;8n+a+
b+c+d).\endalign$$
\endpro
Proof. By Theorem 4.3, $t(a,b,c;k)=N(a,b,c;8k+a+b+c)$ for any
$k\in\Bbb N.$ Now applying Lemma 4.1 (with $m=1$) yields the result.
\par\q
\pro{Theorem 5.4} Suppose $a,c,d,n\in\Bbb N$, $2\nmid a$ and $4\nmid
c$. Then
$$t(a,3a,c,d;n)=2N(4a,12a,c,d;8n+4a+c+d)-2N(4a,12a,c,4d;8n+4a+c+d).$$
\endpro
Proof. By Theorem 3.4, for $m=0,1,2,\ldots$ we have
$t(a,3a,c;m)=2N(4a,12a,c;8m+4a+c).$ Thus,
$$\align t(a,3a,c,d;n)
&=\sum_{w\in\Bbb Z}t(a,3a,c;n-dw(w-1)/2)\\&=2\sum_{w\in\Bbb
Z}N(4a,12a,c;8(n-dw(w-1)/2)+4a+c)
\\&=2\sum_{w\in\Bbb Z}N(4a,12a,c;8n+4a+c+d-d(2w-1)^2)
\\&=2\sum_{w\in\Bbb Z}N(4a,12a,c;8n+4a+c+d-dw^2)
\\&\qq-2\sum_{w\in\Bbb Z}N(4a,12a,c;8n+4a+c+d-d(2w)^2)
\\&=2N(4a,12a,c,d;8n+4a+c+d)-2N(4a,12a,c,4d;8n+4a+c+d)\big).\endalign$$
\pro{Corollary 5.2} Suppose $a,c,d,n\in\Bbb N$, $2\nmid ac$ and $d\e
2,c\mod 4$. Then
$$t(a,3a,c,d;n)=2N(4a,12a,c,d;8n+4a+c+d).$$
\endpro

\pro{Theorem 5.5} Let $a,b,d,n\in\Bbb N$ with $2\nmid ab$. Then
$$t(a,3a,2b,d;n)=\f
23\big(N(a,3a,2b,d;8n+4a+2b+d)-N(a,3a,2b,4d;8n+4a+2b+d)\big).$$
\endpro
Proof. By Theorem 3.1, $t(a,3a,2b;k)=\f 23N(a,3a,2b;8k+4a+2b)$ for
any nonnegative integer $k$. Now the result follows from Lemma 4.1.

\par\q
\pro{Theorem 5.6} Let $a,d,n\in\Bbb N$. Then
$$t(a,3a,9a,d;n)=\f
12\big(N(a,3a,9a,d;8n+13a+d)-N(a,3a,9a,4d;8n+13a+d)\big).$$
\endpro
Proof. By Theorem 2.4, $t(1,3,9;m)=\f 12N(1,3,9;8m+13)$. Thus
applying Lemma 4.1 gives the result.
\par

\par\q

\pro{Theorem 5.7} Let $a,b,c,n\in\Bbb N$ with $2\nmid ab$ and
$n\not\e\f{a+b}2\mod 2$. Then
$$\align &t(a,3a,4b,2c;n)
\\&=\f
23\big(N(a,3a,4b,2c;8n+4a+4b+2c)-N(a,3a,4b,8c;8n+4a+4b+2c)\big).
\endalign$$
\endpro
Proof. By Theorem 3.3, for $m\not\e \f{a+b}2\mod 2$,
$t(a,3a,4b;m)=\f 23N(a,3a,4b;8m+4a+4b).$ Thus,
$$\align &t(a,3a,4b,2c;n)
\\&=\sum_{w\in\Bbb Z}t(a,3a,4b;n-2cw(w-1)/2)=\f 23\sum_{w\in\Bbb
Z}N(a,3a,4b;8(n-cw(w-1))+4a+4b)
\\&=\f 23\sum_{w\in\Bbb Z}N(a,3a,4b;8n+4a+4b+2c-2c(2w-1)^2)
\\&=\f 23\sum_{w\in\Bbb Z}N(a,3a,4b;8n+4a+4b+2c-2cw^2)
\\&\qq-\f 23\sum_{w\in\Bbb Z}N(a,3a,4b;8n+4a+4b+2c-2c(2w)^2)
\\&=\f
23\big(N(a,3a,4b,2c;8n+4a+4b+2c)-N(a,3a,4b,8c;8n+4a+4b+2c)\big).\endalign$$
\par\q
\pro{Theorem 5.8} Let $m,n\in\Bbb N$.
\par $(\t{\rm i})$ If there is a prime divisor $p$
of $2m+1$ such that $\sls{8n+5}p=-1$, then
$$t(1,2,2,4m+2;n)=\f 12N(1,1,4,4m+2;8n+4m+7).$$
\par $(\t{\rm ii})$ If there is a prime divisor $p$
of $2m+1$ such that $\sls{8n+9}p=-1$, then
$$t(1,4,4,4m+2;n)
=\f 14N(1,1,4,4m+2;8n+4m+11).$$
\endpro
Proof. By [S2, Theorem 2.7],
$$t(1,2,2,4m+2;n)=t(1,1,8,8m+4;2n),\q t(1,4,4,4m+2;n)=\f 12t(1,1,8,8m+4;2n+1).$$ Now applying
Theorem 2.3 and [S2, Lemma 2.1] yields the result.

\par\q
\pro{Theorem 5.9} Let $a,b\in\{1,3,5,\ldots\}$. Then for  $n\in\Bbb
N$,
$$t(a,a,2b,4b;n)=N(a,a,b,2b;4n+a+3b)-N(a,a,b,2b;2n+(a+3b)/2).$$
\endpro
Proof. By [S2, Lemma 2.1 and Theorem 2.15],
$$\align t(a,a,2b,4b;n)&=N(a,a,4b,2b;8n+2a+6b)-N(a,a,4b,8b;8n+2a+6b)
\\&=N(a,a,2b,b;4n+a+3b)-N(a,a,2b,4b;4n+a+3b)
\\&=N(a,a,b,2b;4n+a+3b)-N(a,a,b,2b;2n+(a+3b)/2).\endalign$$
This proves the theorem.

\pro{Lemma 5.1} Let $a,b\in\{1,3,5,\ldots\}$ and $n\in\Bbb N$. Then
$$t(a,2a,4a,b;n)=\f 14\Big(N(a,a,a,2b;16n+14a+2b)-N(a,a,2a,b;8n+7a+b)
\Big).$$
\endpro
Proof. By [S2, Theorems 2.6, 2.15 and Lemma 2.1],
$$\align &t(a,2a,4a,b;n)\\&=\f 14t(a,a,4b,2a;4n+3a)
\\&=\f
14\Big(N(a,a,4b,2a;8(4n+3a)+4a+4b)-N(a,a,4b,8a;8(4n+3a)+4a+4b)\Big)
\\&=\f
14\Big(N(a,a,2b,a;4(4n+3a)+2b+2a)-N(a,a,2b,4a;4(4n+3a)+2a+2b)\Big)
\\&=\f
14\Big(N(a,a,2b,a;4(4n+3a)+2a+2b)-N(a,a,b,2a;2(4n+3a)+a+b)
 \Big).\endalign$$
This yields the result.
\par\q
 \pro{Theorem 5.10} Let
$a,b\in\{1,3,5,\ldots\}$ and $n\in\Bbb N$. Then
$$\align &N(a,a,2a,b;2n+a+b)
\\&=\f
13\Big(N(a,a,a,2b;4n+2a+2b)+2N\big
(a,a,a,2b;n+\f{a+b}2\big)\Big)\endalign$$
 and
$$\align t(a,2a,4a,b;n)
=\f 16\Big(N(a,a,a,2b;16n+14a+2b)
-N\big(a,a,a,2b;4n+\f{7a+b}2\big)\Big).
\endalign$$
\endpro
Proof. By [S2, Theorems 2.1, 2.15 and Lemma 2.1],
$$\align &\f 23\Big(N(a,a,a,2b;4n+2a+2b)
-N\big(a,a,a,2b;n+\f{a+b}2\big)\Big)
\\&=t(a,a,2a,4b;n)
\\&=N(a,a,2a,4b;8n+4a+4b)-N(a,a,8a,4b;8n+4a+4b)
\\&=N(a,a,a,2b;4n+2a+2b)-N(a,a,4a,2b;4n+2a+2b)
\\&=N(a,a,a,2b;4n+2a+2b)-N(a,a,2a,b;2n+a+b).
\endalign$$
This yields the first part.
\par Set $n'=4n+3a$. By Lemma 5.1,
$$\align &t(a,2a,4a,b;n)
\\&=\f 14\big(N(a,a,a,2b;16n+14a+2b)-
N(a,a,2a,b;8n+7a+b)\big)
\\&=\f 14\big(N(a,a,a,2b;4n'+2a+2b)-
N(a,a,2a,b;2n'+a+b)\big).
\endalign$$
By the first part,
$$\align&N(a,a,2a,b;2n'+a+b)
\\&=\f
13\Big(N(a,a,a,2b;4n'+2a+2b)+2N\big
(a,a,a,2b;n'+\f{a+b}2\big)\Big).\endalign$$ Hence
$$\align &t(a,2a,4a,b;n)
\\&=\f 14\Big(N(a,a,a,2b;4n'+2a+2b)
\\&\q- \f 13\Big(N(a,a,a,2b;4n'+2a+2b)+2N\big
(a,a,a,2b;n'+\f{a+b}2\big)\Big)\Big).\endalign$$ This yields the
second part. The proof is now complete.
\par\q

\pro{Theorem 5.11} Let $a,b\in\{1,3,5,\ldots\}$ and $n\in\Bbb N$.
Then
$$t(a,a,6a,b;n)=\f 12\Big(N(a,a,3a,2b;16n+16a+2b)-N(a,a,6a,b;
8n+8a+b) \Big).$$
\endpro
Proof. By [S2, Theorems 2.5, 2.15 and Lemma 2.1],
$$\align &t(a,a,6a,b;n)\\&=\f 12t(a,a,6a,4b;4n+3a)
\\&=\f
14\Big(N(a,a,6a,4b;8(4n+3a)+8a+4b)-N(a,a,24a,4b;8(4n+3a)+8a+4b)\Big)
\\&=\f
14\Big(N(a,a,3a,2b;4(4n+3a)+4a+2b)-N(a,a,12a,2b;4(4n+3a)+4a+2b)\Big)
\\&=\f
14\Big(N(a,a,3a,2b;4(4n+3a)+4a+2b)-N(a,a,6a,b;2(4n+3a)+2a+b)
 \Big).\endalign$$
This yields the result.

\par\q
\pro{Theorem 5.12} Suppose $a,b,n\in\Bbb N$, $2\nmid a$ and $b\e
2\mod 4$. Then
$$t(a,a,b,b;n)=N(a,a,b,b;4n+a+b).$$
\endpro
Proof. By (1.3) and (1.7),
$$\aligned &\sum_{n=0}^\infty N(a,a,b,b;n)q^n=\phq a^2\phq b^2\\&
=\big(\phq{4a}^2+4q^{2a}\psq{8a}^2+4q^a\psq {4a}^2\big)\big(\phq
{2b}^2+4q^b\psq{4b}^2\big).
\endaligned$$
Thus, $$\sum_{n=0}^\infty
N(a,a,b,b;4n+2+(-1)^{\f{a-1}2})q^{4n+2+(-1)^{(a-1)/2}} =4q^a\psq
{4a}^2\cdot 4q^b\psq{4b}^2$$ and so
$$ \align &\sum_{n=0}^\infty
N(a,a,b,b;4n+2+(-1)^{\f{a-1}2})q^n\\&=16q^{\f{a+b-2-(-1)^{(a-1)/2}}4}\psq
a^2\psq
b^2=q^{\f{a+b-2-(-1)^{(a-1)/2}}4}\sum_{n=0}^{\infty}t(a,a,b,b;n)q^n.\endalign$$
Hence
$$\align t(a,a,b,b;n)&=N\big(a,a,b,b;4\big(n+
\f{a+b-2-(-1)^{(a-1)/2}}4\big)+2+(-1)^{\f{a-1}2}\big)\\&=N(a,a,b,c;4n+a+b)\endalign$$
as claimed.

\par\q
\pro{Theorem 5.13} Let $a,b\in\{1,3,5,\ldots\}$ and $n\in\Bbb N$.
Then
$$t(a,a,b,b;n)=N(a,a,b,b;4n+a+b)-N(a,a,b,b;2n+(a+b)/2).$$
\endpro
Proof. By (1.3) and (1.7),
$$\aligned &\sum_{n=0}^\infty N(a,a,b,b;n)q^n=\phq a^2\phq b^2\\&
=\big(\phq{4a}^2+4q^{2a}\psq{8a}^2+4q^a\psq {4a}^2\big)\big(\phq
{4b}^2+4q^{2b}\psq{8b}^2+4q^b\psq{4b}^2\big).
\endaligned$$
Thus, if $4\mid a+b$, then
$$\aligned &\sum_{n=0}^\infty
N(a,a,b,b;4n)q^{4n}\\&=\phq{4a}^2\phq{4b}^2+16q^{a+b}\psq{4a}^2\psq{4b}^2+
16q^{2(a+b)}\psq{8a}^2\psq{8b}^2
\endaligned$$ and so
$$\aligned&\sum_{n=0}^\infty N(a,a,b,b;4n)q^n\\&=
\phq a^2\phq b^2+16q^{\f{a+b}4}\psq a^2\psq
b^2+16q^{\f{a+b}2}\psq{2a}^2\psq{2b}^2;
\endaligned$$
if $4\mid a-b$, then
$$\aligned &\sum_{n=0}^\infty N(a,a,b,b;4n+2)q^{4n+2}\\&=
4q^{2a}\psq{8a}^2\phq{4b}^2+4q^{2b}\phq{4a}^2\psq{8b}^2+16q^{a+b}
\psq{4a}^2\psq{4b}^2\endaligned$$ and so
$$\aligned &\sum_{n=0}^\infty N(a,a,b,b;4n+2)q^n\\&=
4q^{\f{a-1}2}\psq{2a}^2\phq{b}^2+4q^{\f{b-1}2}\phq{a}^2\psq{2b}^2+16q^{\f{a
+b-2}4} \psq a^2\psq b^2.\endaligned$$ On the other hand,
$$\aligned&\sum_{n=0}^\infty N(a,a,b,b;n)q^n=\phq a^2\phq b^2\\&
=\big(\phq{2a}^2+4q^a\psq{4a}^2\big)\big(\phq{2b}^2+4q^b\psq{4b}^2
\big)\\&=\phq{2a}^2\phq{2b}^2+16q^{a+b}\psq{4a}^2\psq{4b}^2+4q^a
\psq{4a}^2\phq{2b}^2+4q^b\phq{2a}^2\psq{4b}^2.
\endaligned
$$
Hence,
$$\sum_{n=0}^\infty N(a,a,b,b;2n)q^{2n}
=\phq{2a}^2\phq{2b}^2+16q^{a+b}\psq{4a}^2\psq{4b}^2$$ and
$$\sum_{n=0}^\infty N(a,a,b,b;2n+1)q^{2n+1}
=4q^a\psq{4a}^2\phq{2b}^2+4q^b\phq{2a}^2\psq{4b}^2.$$ Therefore,
$$\sum_{n=0}^\infty N(a,a,b,b;2n)q^n
=\phq a^2\phq b^2+16q^{\f{a+b}2}\psq{2a}^2\psq{2b}^2$$ and
$$\sum_{n=0}^\infty N(a,a,b,b;2n+1)q^n
=4q^{\f{a-1}2}\psq {2a}^2\phq b^2+4q^{\f{b-1}2}\phq a^2\psq{2b}^2.$$
From the above we deduce that for $a+b\e0\pmod 4,$
$$\aligned&\sum_{n=0}^\infty \big(N(a,a,b,b;4n)-N(a,a,b,b;2n)\big)q^n
\\&=16q^{\f{a+b}4}\psq a^2\psq b^2=q^{\f{a+b}4}\sum_{n=0}^\infty
t(a,a,b,b;n)q^n,
\endaligned$$
and for $a\e b\pmod 4,$
$$\aligned&\sum_{n=0}^\infty \big(N(a,a,b,b;4n+2)-N(a,a,b,b;2n+1)\big)q^n
\\&=16q^{\f{a+b-2}4}\psq a^2\psq b^2=q^{\f{a+b-2}4}\sum_{n=0}^\infty
t(a,a,b,b;n)q^n.
\endaligned$$
Now comparing the coefficients of $q^{n+[\f{a+b}4]}$ on both sides
yields the result.

\par\q
\pro{Theorem 5.14} Let $a,b,n\in\Bbb N$, $2\nmid ab$ and $4\mid
a-b$. Then
$$t(a,2a,b,2b;n)=N(a,2a,b,2b;8n+3(a+b))-N(a,2a,b,2b;4n+3(a+b)/2\big).$$
\endpro
Proof. By (1.2)-(1.4),
$$\align &\sum_{n=0}^{\infty}N(a,2a,b,2b;n)q^n
=\phq a\phq{2a}\phq b\phq{2b}
\\&=(\phq{4a}+2q^a\psq{8a})(\phq{8a}+2q^{2a}\psq{16a})
(\phq{4b}+2q^b\psq{8b})(\phq{8b}+2q^{2b}\psq{16b})
\\&=(\phq{4a}\phq{8a}+2q^a\phq{8a}\psq{8a}+2q^{2a}\phq{4a}\psq{16a}+4q^{3a}\psq{8a}\psq{16a})
\\&\qq\times(\phq{4b}\phq{8b}+2q^b\phq{8b}\psq{8b}+2q^{2b}\phq{4b}\psq{16b}+4q^{3b}\psq{8b}\psq{16b}).
\endalign$$
For $a\e b\mod 8$ collecting the terms of the form
$q^{4n+2+(-1)^{(a-1)/2}}$ yields
$$\aligned
&\sum_{n=0}^{\infty}N(a,2a,b,2b;4n+2+(-1)^{(a-1)/2})q^{4n+2+(-1)^{(a-1)/2}}
\\&=4q^{3a}\phq{4b}\phq{8b}\psq{8a}\psq{16a}
+4q^{3b}\phq{4a}\phq{8a}\psq{8b}\psq{16b}
\\&\qq+4q^{a+2b}\phq{8a}\psq{8a}\phq{4b}\psq{16b}
+4q^{2a+b}\phq{4a}\psq{16a}\phq{8b}\psq{8b}.\endaligned\tag 5.1$$
For $a\e 5b\mod 8$ collecting the terms of the form
$q^{4n+2-(-1)^{(a-1)/2}}$ yields
$$\aligned
&\sum_{n=0}^{\infty}N(a,2a,b,2b;4n+2-(-1)^{(a-1)/2})q^{4n+2-(-1)^{(a-1)/2}}
\\&=2q^{a}\phq{8a}\psq{8a}\phq{4b}\phq{8b}+
2q^{b}\phq{8b}\psq{8b}\phq{4a}\phq{8a}
\\&\qq+8q^{3a+2b}\psq{8a}\psq{16a}\phq{4b}\psq{16b}
+8q^{2a+3b}\phq{4a}\psq{16a}\psq{8b}\psq{16b}.\endaligned\tag 5.2$$
On the other hand, using (1.5) we see that
$$\align \phq
a\phq{2a}&=(\phq{16a}+2q^{4a}\psq{32a}+2q^a\psq{8a})(\phq{8a}+2q^{2a}\psq{16a})
\\&=\phq{8a}\phq{16a}+2q^a\phq{8a}\psq{8a}+2q^{2a}\phq{16a}\psq{16a}
\\&\qq+4q^{3a}\psq{8a}\psq{16a}
+2q^{4a}\phq{8a}\psq{32a}+4q^{6a}\psq{16a}\psq{32a}.
\endalign$$
Hence
$$\aligned &\sum_{n=0}^{\infty}N(a,2a,b,2b;n)q^n
=\phq a\phq{2a}\phq b\phq{2b}
\\&=(\phq{8a}\phq{16a}+2q^a\phq{8a}\psq{8a}+2q^{2a}\phq{16a}\psq{16a}
\\&\qq+4q^{3a}\psq{8a}\psq{16a}
+2q^{4a}\phq{8a}\psq{32a}+4q^{6a}\psq{16a}\psq{32a})
\\&\qq\times
(\phq{8b}\phq{16b}+2q^b\phq{8b}\psq{8b}+2q^{2b}\phq{16b}\psq{16b}
\\&\qq+4q^{3b}\psq{8b}\psq{16b}
+2q^{4b}\phq{8b}\psq{32b}+4q^{6b}\psq{16b}\psq{32b})
.\endaligned\tag 5.3$$ For $a\e b\mod 8$ collecting the terms of the
form $q^{8n+4+2(-1)^{(a-1)/2}}$ in (5.3) yields
$$\align &\sum_{n=0}^{\infty}N(a,2a,b,2b;8n+4+2(-1)^{(a-1)/2})q^{8n+4+2(-1)^{(a-1)/2}}
\\&=4q^{6a}\psq{16a}\psq{32a}\phq{8b}\phq{16b}+4q^{6b}\phq{8a}\phq{16a}\psq{16b}\psq{32b}
\\&\qq+4q^{2a+4b}\phq{16a}\psq{16a}\phq{8b}\psq{32b}+4q^{4a+2b}\phq{8a}\psq{32a}
\phq{16b}\psq{16b}
\\&\qq+16q^{3a+3b}\psq{8a}\psq{16a}\psq{8b}\psq{16b}.
\endalign$$
and so
$$\align
&\sum_{n=0}^{\infty}N(a,2a,b,2b;8n+4+2(-1)^{(a-1)/2})q^{4n+2+(-1)^{(a-1)/2}}
\\&=4q^{3a}\psq{8a}\psq{16a}\phq{4b}\phq{8b}+4q^{3b}\phq{4a}\phq{8a}\psq{8b}\psq{16b}
\\&\qq+4q^{a+2b}\phq{8a}\psq{8a}\phq{4b}\psq{16b}+4q^{2a+b}\phq{4a}\psq{16a}
\phq{8b}\psq{8b}
\\&\qq+16q^{3(a+b)/2}\psq{4a}\psq{8a}\psq{4b}\psq{8b}.
\endalign$$
This together with (5.1) yields
$$\aligned &\sum_{n=0}^{\infty}(N(a,2a,b,2b;8n+4+2(-1)^{(a-1)/2})-N(a,2a,b,2b;4n+2+(-1)^{(a-1)/2})q^{4n+2+(-1)^{(a-1)/2}}
\\&=16q^{3(a+b)/2}\psq{4a}\psq{8a}\psq{4b}\psq{8b}
\endaligned\tag 5.4$$
Since $\f{3(a+b)}2\e 2+(-1)^{(a-1)/2}\mod 4$, substituting $q$ with
$q^{\f 14}$ we get
$$\aligned &\sum_{n=0}^{\infty}(N(a,2a,b,2b;8n+3(a+b))-N(a,2a,b,2b;4n+3(a+b)/2))q^n
\\&=16\psq{a}\psq{2a}\psq{b}\psq{2b}=\sum_{n=0}^{\infty}t(a,2a,b,2b;n)q^n,
\endaligned\tag 5.5$$
which yields the result in this case.
\par For $a\e 5b\mod 8$ collecting the terms of the
form $q^{8n+4-2(-1)^{(a-1)/2}}$ in (5.3) yields
$$\align &\sum_{n=0}^{\infty}N(a,2a,b,2b;8n+4-2(-1)^{(a-1)/2})q^{8n+4-2(-1)^{(a-1)/2}}
\\&=2q^{2a}\phq{16a}\psq{16a}\phq{8b}\phq{16b}+2q^{2b}\phq{8a}\phq{16a}\phq{16b}\psq{16b}
\\&\qq+8q^{4a+6b}\phq{8a}\psq{32a}
\psq{16b}\psq{32b}+8q^{6a+4b}\psq{16a}\psq{32a}\phq{8b}\psq{32b}
\\&\qq+16q^{3a+3b}\psq{8a}\psq{16a}\psq{8b}\psq{16b}.
\endalign$$
and so
$$\align &\sum_{n=0}^{\infty}N(a,2a,b,2b;8n+4-2(-1)^{(a-1)/2})q^{4n+2-(-1)^{(a-1)/2}}
\\&=2q^{a}\phq{8a}\psq{8a}\phq{4b}\phq{8b}+2q^{b}\phq{4a}\phq{8a}\phq{8b}\psq{8b}
\\&\qq+8q^{2a+3b}\phq{4a}\psq{16a}
\psq{8b}\psq{16b}+8q^{3a+2b}\psq{8a}\psq{16a}\phq{4b}\psq{16b}
\\&\qq+16q^{3(a+b)/2}\psq{4a}\psq{8a}\psq{4b}\psq{8b}.
\endalign$$
This together with (5.2) yields
$$\aligned &\sum_{n=0}^{\infty}(N(a,2a,b,2b;8n+4-2(-1)^{(a-1)/2})-N(a,2a,b,2b;4n+2-(-1)^{(a-1)/2})q^{4n+2-
(-1)^{(a-1)/2}}
\\&=16q^{3(a+b)/2}\psq{4a}\psq{8a}\psq{4b}\psq{8b}
\endaligned\tag 5.6$$
Since $\f{3(a+b)}2\e 2-(-1)^{(a-1)/2}\mod 4$, substituting $q$ with
$q^{\f 14}$ in (5.6) yields (5.5). Hence the result is true when
$a\e 5b\mod 8$. The proof is now complete.

\par\q
\pro{Theorem 5.15} Let $n\in\Bbb N$. Then
$$t(1,1,1,6;n)=\f 16\Big(N(1,1,1,6;32n+36)-N(1,1,1,6;8n+9)\Big).
$$\endpro
Proof. By Theorem 5.11,
$$t(1,1,1,6;n)=\f 12\Big(N(1,1,2,3;16n+18)-N(1,1,1,6;8n+9)\Big).$$
By Theorem 5.10,
$$N(1,1,2,3;2m+4)=\f 13\Big(N(1,1,1,6;4m+8)+2N(1,1,1,6;m+2)\Big).$$
Thus,
$$\align N(1,1,2,3;16n+18)&=N(1,1,2,3;2(8n+7)+4)
\\&=\f 13\Big(N(1,1,1,6;32n+36)+2N(1,1,1,6;8n+9)
\Big).\endalign$$
 Now combining
all the above gives the result.
\par\q
\pro{Lemma 5.2} For
$|q|<1$ we have
$$\align \varphi(q)^3=&\phq{16}^3+6q^4\phq{16}^2\psq{32}+12q^8\phq{16}\psq{32}^2+8q^{12}\psq{32}^3
\\&+6q\phq 8^2\psq 8+24q^5\psq 8\psq{16}^2+12q^2\phq{16}\psq
8^2\\&+24q^6\psq 8^2\psq{32}+8q^3\psq 8^3.\endalign$$
\endpro
Proof. By [S2, Lemma 2.2],
 $$\varphi(q)^3=\phq 4^3+6q\phq 4\psq 4^2+12q^2\psq 4^2\psq 8+8q^3\psq
 8^3.$$
 By (1.1)-(1.2),
 $$\align &\phq 4=\phq{16}+2q^4\psq{32},\\& \psq 4^2=\phq 4\psq
 8=(\phq{16}+2q^4\psq{32})\psq 8,
 \\&\phq 4\psq 4^2=\phq 4^2\psq 8=(\phq 8^2+4q^4\psq{16}^2)\psq 8.\endalign$$
 Thus,
$$\align \varphi(q)^3&=(\phq{16}+2q^4\psq{32})^3+6q(\phq 8^2+4q^4\psq{16}^2)\psq 8
\\&\q+12q^2(\phq{16}+2q^4\psq{32})\psq 8^2+8q^3\psq
 8^3.\endalign$$
 This yields the result.
\par\q
\pro{Lemma 5.3} For $|q|<1$ we have
$$\varphi(q)\phq 7-\phq 2\phq{14}=2q\psi(q)\psq 7-4q^2\psq
2\psq{14}+4q^4\psq 4\psq{28}.$$
\endpro
Proof. Let $\sls am$ be the Legendre-Jacobi-Kronecker symbol, and
let $[q^n]f(q)$ be the coefficient of $q^n$ in the power series
expansion of $f(q)$. Suppose $n\in\Bbb N$ and $n=2^{\alpha}n_0\
(2\nmid n_0)$. Set
$$R(n)=|\{(x,y)\in\Bbb Z\times \Bbb Z\bigm|n=x^2+7y^2\}\qtq{and}T(n)=\sum_{k\mid n,2\nmid k}\Ls {-7}k
=\sum_{k\mid n_0}\Ls{-7}k.$$
 By [Be1, pp.302-303],
$$q\psi(q)\psq
7=\sum_{n=1}^{\infty}\Ls{-28}n\f{q^n}{1-q^n}=
\sum_{n=1}^{\infty}\Ls{-28}n\sum_{k=1}^{\infty}q^{kn}
=\sum_{n=1}^{\infty}T(n)q^n.$$ By [SW, Theorem 4.1],
$$R(n)=\cases 2\sum_{k\mid n}\sls{-7}k&\t{if $2\nmid n$,}
\\0&\t{if $4\mid n-2$,}
\\ 2\sum_{k\mid \f n4}\sls{-7}k&\t{if $4\mid n$.}
\endcases$$
If $2\nmid n$, then
$$\align &[q^n](\varphi(q)\phq 7-\phq 2\phq{14})
\\&=[q^n]\varphi(q)\phq 7=R(n)=2\sum_{k\mid n}\Ls{-7}k
=2T(n)=[q^n](q\psi(q)\psq 7)\\&=[q^n](2q\psi(q)\psq 7-4q^2\psq
2\psq{14}+4q^4\psq 4\psq{28}).\endalign$$ If $4\mid n-2$, then
$$\align &[q^n](\varphi(q)\phq 7-\phq 2\phq{14})\\&=R(n)-R\Big(\f n2\Big)=-R\Big(\f
n2\Big)=-2\sum_{k\mid \f n2}\Ls{-7}k
\\&=2\sum_{k\mid n,2\nmid k}\Ls{-7}k-4\sum_{k\mid \f n2}\Ls{-7}k
=2T(n)-4T\ls n2\\&=[q^n](2q\psi(q)\psq 7-4q^2\psq 2\psq{14}+4q^4\psq
4\psq{28}).\endalign$$ If $4\mid n$, then
$$\align &[q^n](\varphi(q)\phq 7-\phq 2\phq{14})=R(n)-R\Big(\f n2\Big)
\\&=\cases 2\sum_{k\mid \f n4}\sls{-7}k-2\sum_{k\mid\f n8}\sls{-7}k
=2\sum_{k\mid n_0}\sls{-7}{2^{\alpha-2}k}=2\sum_{k\mid
n_0}\sls{-7}k&\t{if $8\mid n$,}
\\2\sum_{k\mid \f n4}\sls{-7}k&\t{if $8\mid n-4$}
\endcases
\\&=2\sum_{k\mid n_0}\Ls{-7}k
=2T(n)=2T(n)-4T\Ls n2+4T\Ls n4
\\&=[q^n](2q\psi(q)\psq 7-4q^2\psq 2\psq{14}+4q^4\psq
4\psq{28}).
\endalign$$
Summarizing the above proves the lemma.
\par\q
 \pro{Theorem 5.16} Let $n\in\Bbb N$. Then
$$t(1,1,1,7;n)=4N(1,1,1,7;4n+5)-2N(1,1,1,7;8n+10)$$
and
$$t(1,7,7,7;n)=4N(1,7,7,7;4n+11)-2N(1,7,7,7;8n+22).$$
\endpro
Proof. By Lemma 5.2,
$$\align&\sum_{n=0}^{\infty}N(1,1,1,7;n)q^n=\varphi(q)^3\phq 7
\\&=(\phq{16}^3+6q^4\phq{16}^2\psq{32}+12q^8\phq{16}\psq{32}^2+8q^{12}\psq{32}^3
\\&+6q\phq 8^2\psq 8+24q^5\psq 8\psq{16}^2+12q^2\phq{16}\psq
8^2\\&+24q^6\psq 8^2\psq{32}+8q^3\psq
8^3)\\&\q\times(\phq{112}+2q^7\psq{56}+2q^{28}\psq{224}).\endalign$$
 Thus,
$$\align &\sum_{n=0}^{\infty}N(1,1,1,7;8n+2)q^{8n+2}
\\&=12q^2\phq{16}\psq 8^2\phq{112}+16q^{10}\psq 8^3\psq{56}+48q^{34}\psq 8^2\psq{32}\psq{224}
\endalign$$
and so
$$\align &\sum_{n=0}^{\infty}N(1,1,1,7;8n+2)q^n
\\&=12\phq{2}\psi(q)^2\phq{14}+16q\psi(q)^3\psq{7}+48q^{4}\psi(q)^2\psq{4}\psq{28}.
\endalign$$
On the other hand, using [S2, Lemma 2.2] we see that
$$\align&\sum_{n=0}^{\infty}N(1,1,1,7;n)q^n=\varphi(q)^3\phq 7
\\&=(\phq 4^3+6q\phq 4\psq 4^2+12q^2\psq 4^2\psq 8+8q^3\psq
 8^3)(\phq {28}+2q^7\psq{56}).\endalign$$
Thus,
$$\sum_{n=0}^{\infty}N(1,1,1,7;4n+1)q^{4n+1}
=6q\phq 4\psq 4^2\phq{28} +24q^9\psq 4^2\psq 8\psq{56}$$ and so
$$\sum_{n=0}^{\infty}N(1,1,1,7;4n+1)q^n
=6\varphi(q)\psi(q)^2\phq{7} +24q^2\psi(q)^2\psq 2\psq{14}.$$ Hence
applying Lemma 5.3 gives
$$\align&\sum_{n=0}^{\infty}(2N(1,1,1,7;4n+1)-N(1,1,1,7;8n+2))q^n
\\&=\psi(q)^2(12\varphi(q)\phq 7-12\phq 2\phq {14}+48q^2\psq 2\psq{14}-48q^4\psq
4\psq {28}\\&\q-16q\psi(q)\psq 7) \\&=8q\psi(q)^3\psq 7=\f
12\sum_{n=0}^{\infty}t(1,1,1,7;n)q^{n+1}.\endalign$$ Therefore
$t(1,1,1,7;n)=4N(1,1,1,7;4n+5)-2N(1,1,1,7;8n+10)$. The remaining
part of the theorem can be proved similarly.
\par\q
\pro{Theorem 5.17} Let $n\in\Bbb N$. Then
$$t(1,2,6,6;n)=2N(1,2,6,6;8n+15)-N(1,2,6,6;16n+30)$$
and
$$t(2,2,3,6;n)=2N(2,2,3,6;8n+13)-N(2,2,3,6;16n+26).$$
\endpro
Proof. By (1.3) and (1.5),
$$\align &\sum_{n=0}^{\infty}N(1,2,6,6;n)q^n=\varphi(q)\phq 2\phq
6^2
\\&=(\phq{16}+2q^4\psq{32}+2q\psq 8)(\phq
8+2q^2\psq{16})\\&\q\times(\phq{24}^2+4q^{12}\psq{48}^2+4q^6\psq{24}^2)
\\&=(\phq 8\phq{16}+2q\phq 8\psq 8+2q^2\phq{16}\psq{16}+4q^3\psq
8\psq{16}+2q^4\phq 8\psq{32}\\&\q+4q^6\psq
{16}\psq{32})\times(\phq{24}^2+4q^{12}\psq{48}^2+4q^6\psq{24}^2) .
\endalign$$
Collecting the terms of the form $q^{8n+7}$ yields
$$\sum_{n=0}^{\infty}N(1,2,6,6;8n+7)q^{8n+7}=16q^{15}\psq
8\psq{16}\psq{48}^2+8q^7\phq 8\psq 8\psq{24}^2$$ and so
$$\sum_{n=0}^{\infty}N(1,2,6,6;8n+7)q^n=16q\psi(q)\psq 2\psq
6^2+8\varphi(q)\psi(q)\psq 3^2.$$
\par If $16n+30=x^2+2y^2+6z^2+6w^2$ for $x,y,z,w\in\Bbb Z$, then
$2\mid x$ and so $16n+30=4x^2+2y^2+6z^2+6w^2$ for some
$x,y,z,w\in\Bbb Z$. That is, $8n+15=2x^2+y^2+3z^2+3w^2$ for some
$x,y,z,w\in\Bbb Z$. Hence
$$N(1,2,6,6;16n+30)=N(1,2,3,3;8n+15).$$
By (1.5),
$$\align &\sum_{n=0}^{\infty}N(1,2,3,3;n)q^n=\varphi(q)\phq 2\phq
3^2
\\&=(\phq{16}+2q^4\psq{32}+2q\psq 8)(\phq
8+2q^2\psq{16})(\phq{48}+2q^{12}\psq{96}+2q^3\psq{24})^2
\\&=(\phq 8\phq{16}+2q\phq 8\psq 8+2q^2\phq{16}\psq{16}+4q^3\psq
8\psq{16}+2q^4\phq 8\psq{32}\\&\q+4q^6\psq {16}\psq{32})\times
(\phq{48}^2+4q^{24}\psq{96}^2+4q^6\psq{24}^2\\&\q+4q^{12}\phq{48}\psq{96}+4q^3\phq{48}\psq{24}+8q^{15}
\psq{24}\psq{96}).
\endalign$$
Collecting the terms of the form $q^{8n+7}$ yields
$$\align&\sum_{n=0}^{\infty}N(1,2,3,3;8n+7)q^{8n+7}
\\&=8q^7\phq 8\psq 8\psq{24}^2+16q^{15}\psq
8\psq{16}\phq{48}\psq{96}
\\&\q+8q^7\phq 8\psq{32}\phq{48}\psq{24}
+8q^{15}\phq 8\phq {16}\psq{24}\psq{96}\endalign$$ and so
$$\align&\sum_{n=0}^{\infty}N(1,2,3,3;8n+7)q^n
\\&=8\varphi(q)\psi(q)\psq 3^2+16q\psi(q)\psq 2\phq 6\psq{12}
\\&\q+8\varphi(q)\psq 4\phq 6\psq 3+8q\varphi(q)\phq 2\psq
3\psq{12}.\endalign$$ Applying (1.1) and (1.4) we see that
$$\align&\sum_{n=0}^{\infty}N(1,2,3,3;8n+7)q^n
\\&=8\varphi(q)\psi(q)\psq 3^2+16q\psi(q)\psq 2\psq
6^2+8\varphi(q)\psq 3\psi(q)\psq 3 \\&=16\varphi(q)\psi(q)\psq
3^2+16q\psi(q)\psq 2\psq 6^2.
\endalign$$
Hence
$$\align&\sum_{n=0}^{\infty}(2N(1,2,6,6;8n+7)-N(1,2,3,3;8n+7))q^n
\\&=32q\psi(q)\psq 2\psq
6^2+16\varphi(q)\psi(q)\psq 3^2 -16\varphi(q)\psi(q)\psq
3^2-16q\psi(q)\psq 2\psq 6^2\\&=16q\psi(q)\psq 2\psq 6^2
=\sum_{n=0}^{\infty}t(1,2,6,6;n)q^{n+1}.
\endalign$$
Comparing the coefficients of $q^{n+1}$ on both sides yields
$$\align t(1,2,6,6;n)&=2N(1,2,6,6;8n+15)-N(1,2,3,3;8n+15)\\&=2N(1,2,6,6;8n+15)-N(1,2,6,6;16n+30).
\endalign$$ The remaining part can be proved similarly.

\section*{6. Some open conjectures}
\par\q
\pro{Conjecture 6.1} Let $n\in\Bbb N$. For
$(a,b,c)=(1,1,7),(1,1,15),
 (1,7,7),(1,7,15),$ $
 (1,9,15),(1,15,15),(1,15,25)$
we have
$$t(a,b,c;n)=\f 12\Big(N(a,b,c;4(8n+a+b+c))
-N(a,b,c;8n+a+b+c)\Big).$$
\endpro

\pro{Conjecture 6.2} Let $n\in\Bbb N$. For
$(a,b,c)=(1,3,5),(1,3,7),(1,3,15),$ $(1,3,21),$
$(1,5,15),(1,7,21),(3,5,9),(3,5,15),(3,7,21)$ we have
$$t(a,b,c;n)=\f 12\Big(3N(a,b,c;8n+a+b+c)
-N(a,b,c;4(8n+a+b+c))\Big).$$
\endpro
\pro{Conjecture 6.3} Let $n\in\Bbb N$ with
 $2\mid n$. Then
 $$\align &t(1,2,15;n)=\f
 12\big(N(1,2,15;4(8n+18))
 -N(1,2,15;8n+18)\big),
\\&t(1,15,18;n)=\f
 12\big(N(1,15,18;4(8n+34))
 -N(1,15,18;8n+34)\big),
 \\&t(1,15,30;n)=\f
 12\big(N(1,15,30;4(8n+46))
 -N(1,15,30;8n+46)\big).
 \endalign$$

\pro{Conjecture 6.4} Let $n\in\Bbb N$ with
 $n\e 0,2\mod 3$. Then

$$t(1,1,27;n)=\f 12\big(N(1,1,27;4(8n+29))-N(1,1,27;8n+29)\big)$$.
\endpro

\pro{Conjecture 6.5} Let $n\in\Bbb N$. If
$(a,b,c,d)=(1,1,6,9),(1,3,3,6),(1,6,9,9),(2,3,3,3)$, then
$$t(a,b,c,d;n)=\f
16\big(N(a,b,c,d;4(8n+a+b+c+d))-N(a,b,c,d;8n+a+b+c+d)\big).$$
\endpro
\par By calculations with Maple, the relation $$t(a,b,c,d;n)=\f
16\big(N(a,b,c,d;4(8n+a+b+c+d))-N(a,b,c,d;8n+a+b+c+d)\big)$$ holds
for
$$\align
(a,b,c,d)=&(1,1,1,2),(1,1,1,3),(1,1,1,4),(1,1,1,5),(1,1,1,6),(1,1,2,2),(1,1,2,3),
\\&(1,1,2,4),(1,1,3,3),(1,1,3,9),(1,1,6,9),(1,2,2,2),(1,2,2,3),(1,3,3,3),
\\&(1,3,3,6),(1,3,6,6),(1,3,9,9),(1,6,9,9),(2,3,3,3).\endalign$$

\pro{Conjecture 6.6} Let $n\in\Bbb N$. Then
$$\align t(1,1,1,7;n)&=\f 13\Big(N(1,1,1,7;16n+20)-N(1,1,1,7;4n+5)\Big)
\\&=\f 27\Big(N(1,1,1,7;32n+40)-2N(1,1,1,7;8n+10)\Big)\endalign$$
and
$$\align t(1,7,7,7;n)&
=\f 13\big(N(1,7,7,7;16n+44)-N(1,7,7,7;4n+11)\big)
\\&=\f 27\big(N(1,7,7,7;32n+88)-2N(1,7,7,7;8n+22)\big).\endalign$$
\endpro

\pro{Conjecture 6.7} Let $n\in\Bbb N$. Then $n$ is represented by
$\f{x(x-1)}2+\f{y(y-1)}2+6\f{z(z-1)}2$ if and only if $n\not\e
2\cdot 3^{2r-1}-1\mod{3^{2r}}$ for $r=1,2,3,\ldots$.
\endpro
\par\q
\par Suppose $m=x^2+y^2+6z^2$ for $x,y,z\in\Bbb Z$
and $m=3^{2r-1}m_0$ with $r\ge 2$ and $3\nmid m_0$.
 Since $x^2+
y^2\e x^2+y^2+6z^2=m\e 0\mod 3$ we see that
 $3\mid x$
and $3\mid y$. Hence $m=(3x)^2+(3y)^2+6z^2$ for some $x,y,z\in \Bbb
Z$. Since $9\mid m$, we must have $3\mid z$ and so
$m=(3x)^2+(3y)^2+6(3z)^2$ for some $x,y,z\in \Bbb Z$. That is, $\f
m9=x^2+y^2+6z^2$ for some $x,y,z\in \Bbb Z$. Repeating the procure,
we derive that $3m_0=\f{m}{3^{2r-2}}=x^2+y^2+6z^2$ for some
$x,y,z\in \Bbb Z$. Hence $3\mid x^2+y^2$. This yields $3\mid x$ and
$3\mid y$. Therefore, $3m_0\e 6z^2\mod 9$ and so $m_0 \e 2z^2 \mod
3$. Since $3\nmid m_0$ we must have $3\nmid z$ and so $z^2\e 1\mod
3$. Hence $m_0\e 2\mod 3$ and $m=3^{2r-1}m_0 \e 2\cdot
3^{2r-1}\mod{3^{2r}}$. This shows that for $m\e 3^{2r-1}
\mod{3^{2r}}$, $m$ is not represented by $x^2+y^2+6z^2$. If
$n=\f{x(x-1)}2+\f{y(y-1)}2+6\f{z(z-1)}2$ for $x,y,z\in \Bbb Z$, then
$8n+8=x^2+y^2+6z^2$ with odd integers $x,y$ and $z$. Hence $8(n+1)
\not\e 3^{2r-1}\mod{3^{2r}}$ for every positive integer $r$. This
yields $n\not\e 2\cdot 3^{2r-1}-1\mod{3^{2r}}$ for $r=1,2,3,
\ldots$.

\end{document}